\documentclass[reqno]{amsart}
\usepackage{amsmath,amsthm,amssymb}

\theoremstyle{plain}
\newtheorem{thm}{Theorem}[section]
\newtheorem{prop}[thm]{Proposition}

\newtheorem{cor}[thm]{Corollary}
\theoremstyle{definition}
\newtheorem{rem}{Remark}
\newtheorem{defn}[thm]{Definition}
\newtheorem{eg}[thm]{Example}
\numberwithin{equation}{section}

\newcommand{\bthm}{\begin{thm}}
\newcommand{\ethm}{\end{thm}}
\newcommand{\bprop}{\begin{prop}}
\newcommand{\eprop}{\end{prop}}
\newcommand{\bcor}{\begin{cor}}
\newcommand{\ecor}{\end{cor}}
\newcommand{\bca}{\begin{cases}}
\newcommand{\eca}{\end{cases}}
\newcommand{\brem}{\begin{rem}}
\newcommand{\erem}{\end{rem}}
\newcommand{\bpm}{\begin{pmatrix}}
\newcommand{\epm}{\end{pmatrix}}
\newcommand{\bdefn}{\begin{defn}}
\newcommand{\edefn}{\end{defn}}
\newcommand{\bsub}{\begin{subtitle}}
\newcommand{\esub}{\end{subtitle}}
\newcommand{\ben}{\begin{enumerate}}
\newcommand{\een}{\end{enumerate}}
\newcommand{\beg}{\begin{eg}}
\newcommand{\eeg}{\end{eg}}

\def\bs{\bigskip}
\def\ni{\noindent}

\def\p{\partial}

\def\diag{{\rm diag\/}}
\def\ti{\tilde}

\def \a {\alpha}
\def \b {\beta}
\def \d {\delta}

\def \e {\epsilon}
\def \g {\gamma}

\def \l {\lambda}

\def \n {\,\vert\,}

\def \o {\theta}

\def\W{\Omega}
\def \s {\sigma}
\def\bu{\bullet}

\def\R{\mathbb{R} }
\def\C{\mathbb{C}}

\def\cg{{\mathcal {G}}}

\def\ck{{\mathcal {K}}}
\def\cm{{\mathcal {M}}}

\def\co{{\mathcal {O}}}
\def\cp{{\mathcal {P}}}

\def\cu{{\mathcal {U}}}

\def \li{\langle}
\def \ri{\rangle}
\def \n {\ \vert\ }
\def\tr{{\rm tr}}
\def\sech{{\rm sech\ }}

\begin{document}


\author{ Chuu-Lian Terng$^1$}\thanks{$^1$Research supported
in  part by
NSF Grant DMS-0306446}
\address{Department of Mathematics\\
Northeastern University\\Boston, MA 02115}
\email{terng@neu.eud}
\author{Karen Uhlenbeck$^2$}
\thanks{$^2$Research supported in part
by  Sid Richardson
Regents' Chair Funds, University of Texas system and NSF Grant DMS-0305505}
\address{Department of Mathematics\\
University of Texas at Austin\\
Austin, TX 78712}
\email{uhlen@math.utexas.edu}

\title{$1+1$ wave maps into symmetric spaces}

\begin{abstract}

We explain how to apply techniques from integrable systems to construct $2k$-soliton
homoclinic wave  maps from the periodic Minkowski space $S^1\times R^1$ to a compact
Lie group, and more generally to a compact symmetric space.  We give a correspondence
between solutions of the $-1$ flow equation associated to a compact Lie group $G$ and
wave maps into $G$.  We use B\"acklund transformations to construct explicit
$2k$-soliton breather solutions for the $-1$ flow equation and show that the
corresponding wave maps are periodic and homoclinic. The compact symmetric space
$G/K$ can be embedded as a totally geodesic submanifold of $G$ via the Cartan
embedding.  We prescribe the constraint condition for the $-1$ flow equation associated
to $G$ which insures that the corresponding wave map into $G$ actually lies in $G/K$.    
For example, when
$G/K=SU(2)/SO(2)=S^2$, the constrained $-1$-flow equation associated to $SU(2)$ has the
sine-Gordon equation (SGE) as a subequation and classical breather solutions of the SGE are
$2$-soliton breathers.  Thus our result generalizes the result of Shatah and Strauss that a classical
breather solution of the SGE gives rise to a periodic homoclinic wave map  to $S^2$. When the group
$G$ is non-compact, the bi-invariant metric on $G$ is pseudo-Riemannian and B\"acklund
transformations of a smooth solution often are singular.  We use B\"acklund transformations to
show that there exist smooth initial data with constant boundary conditions and finite energy such
that the Cauchy problem for wave maps from $R^{1,1}$ to the pseudo-Riemannian manifold
$SL(2,R)$ develops singularities in finite time.

\end{abstract}


\maketitle

\section{Introduction}\label{fa}

A smooth map $s:M\to N$ between two pseudo-Riemannian manifolds is called {\it harmonic\/}
 if it is a critical point of the  functional 
\begin{equation}\label{ev}
J(s)=\frac{1}{2}\int_M \li ds_p, ds_p\ri_p dv,
\end{equation}
where $\li \ , \ri_p$ is the natural bilinear form induced from the metrics on $TM_p$ and
$TN_{s(p)}$, and $dv$ is the volume form of $M$ associated to its metric. When the
domain manifold is Riemannian, the Euler-Lagrange equation of
$J$ is elliptic, and is the natural non-linear generalization of the Laplace-Beltrami equation.   There is an
extensive literature in geometry and physics on elliptic harmonic maps.  When the domain manifold
 is the Lorentz  space $M=\R^{n,1}$, the equations are non-linear wave equations,
and the solution maps are referred to as {\it wave maps\/}. When the target
manifold $N$ is a Lie Group $G$, these equations have a particularly simple
form. For  $M= \R^{1,1}$, the equation for $s: \R^{1,1} \to G$ is
\begin{equation} \label{ae}
(s^{-1}s_t)_t= (s^{-1}s_x)_x,
\end{equation}
and a solution is called a $1+1$ wave map into $G$. 
This equation in light cone (characteristic)
coordinates $$\xi=\frac{x+t}{2}, \ \   \eta=\frac{x-t}{2},$$  takes the form
\begin{equation} \label{am}
(s^{-1}s_\xi)_\eta = -(s^{-1}s_\eta)_\xi,
\end{equation}
which can be encoded
in a Lax pair, i.e., $s$ is a $1+1$ wave map if and only if 
$$\left[\frac{\p}{\p \xi}+ \frac{(1-\l)}{2}\ s^{-1} s_\xi, \ \ \frac{\p}{\p \eta} +
\frac{(1-\l^{-1})}{2}\ s^{-1}s_\eta\right]=0$$
for all $\l\in \C\setminus \{0\}$.  Because of this Lax formulation, the $1+1$  wave map equation
\eqref{ae} is an {\it integrable system\/}. 

  In this paper we describe how to
apply methods from integrable systems to construct periodic and homoclinic wave maps to  Lie groups and
more generally to symmetric spaces.  We will use a closely related integrable non-linear wave equation, the
$-1$ flow equation associated to $G$.  This is the first order semi-linear wave system for $(a,
u,v):\R^{1,1}\to
\cg\times
\cg\times \cg$:
\begin{equation}\label{eu}
\bca a_t= a_x,&\cr
u_t=u_x- [a,v],&\cr
v_t= -v_x -[u,v],
\eca 
\end{equation}
where $\cg$ is the Lie algebra of $G$.  The $-1$ flow equation also has a Lax pair, namely
$(a,u,v)$ is a solution of the $-1$ flow equation if and only if 
$$\left[\frac{\p}{\p \xi} + a\l + u, \ \  \frac{\p}{\p \eta} + \l^{-1} v  \right]=0$$
for all $\l\in \C\setminus \{0\}$, where $\xi, \eta$ are characteristic coordinates. 
The name $-1$ flow comes from the standard convention in the theory of soliton equations and we give a brief explanation next. 

There is a hierarchy of soliton flows associated to each  Lie group $G$.  The Lax
pair of the $j$-th flow in the $G$-hierarchy is of the form
$$\left[\frac{\p}{\p x} + a\l + u, \ \ \frac{\p}{\p t} + b\l^j + Q_1\l^{j-1} +\cdots +
Q_j\right]=0$$
with $a, b, u, Q_i$ in $\cg$.  
For example,  the focusing non-linear Schr\"odinger equation (NLS) 
$$q_t=\frac{i}{2}(q_{xx}+ 2|q|^2 q)$$ and the complex modified KdV equation $$q_t= -\frac{1}{4}(q_{xxx} + 6|q|^2 q_x)$$  are the
second and third flows in the $SU(2)$-hierarchy.  

The Lax pairs of the wave map equation and the $-1$ flow equation are gauge equivalent, which in
turn gives an equivalence between solutions of the $-1$ flow  equation associated to $G$ and wave maps
$s: \R^{1,1} \to G$ with $s(0,0)={\rm I}$ the identity. 
 This is analogous to the Hasimoto
transformation between the focusing NLS and the Heisenberg magnetic model
equation.   If a solution of the
$-1$ flow equation associated to $G$ satisfies a certain constraint (a {\it reality condition\/}) coming from an
involution $\sigma$ of
$G$, then the corresponding wave map to $G$ is in fact a wave map to the symmetric space $G/K$ (here $K$ is
the fixed point set of $\sigma$).   When
$G=SU(2)$ and
$\sigma(g)=(g^t)^{-1}$, the constrained
$-1$-flow equation is equivalent to the equation for wave maps from
$\R^{1.1}$ to $S^2$.  Moreover, this constrained $-1$ flow equation associated to $SU(2)$
contains the sine-Gordon equation (SGE),
$$q_{tt}-q_{xx}= \sin q,$$ as a subequation.

Shatah and Strauss prove in \cite{ShaStr96} that wave maps into $S^2$ corresponding 
to classical breather solutions of the SGE  are
homoclinic, when viewed as wave maps from $S^1\times \R^1$ to $S^2$, in the sense that their limits as
$t\to \infty$ and as $t\to -\infty$ are the same.  Our study of periodic $1+1$ wave maps into
symmetric spaces was inspired by their paper.  

 To explain the method we use to construct homoclinic wave maps from $S^1\times \R$ to a 
symmetric space, we need to give a brief review of B\"acklund transformations in integrable systems (cf.
\cite{TerUhl00a}). First note that given smooth $gl(n,\C)$ valued maps $A, B$ on $\R^2$, the condition
that $A, B$ satisfy
 $$\left[\frac{\p}{\p \xi} + A,\ \  \frac{\p}{\p \eta}+ B\right]=0$$ is equivalent to the existence 
of the trivialization $E$ such that 
 $$E_\xi= EA, \quad E_\eta= EB, \quad E(0,0)={\rm I\/}.$$  Now let $(a, u,v)$ be a solution of the 
$-1$ flow, and $E(\xi, \eta, \l)$ the trivialization of the corresponding Lax pair, i.e., 
 $$E_\xi= E\ (a\l + u), \quad E_\eta= \l^{-1} E v, \quad E(0,0,\l)={\rm I\/}.$$
 Since the coefficients of the above differential equation are holomorphic in the parameter 
$\l\in \C\setminus \{0\}$, $E(\xi, \eta,\l)$ is holomorphic in $\l\in \C\setminus \{0\}$.  
The basic idea of
a B\"acklund transformation is that given a linear fractional map from $S^2=\C\cup \{\infty\}$ to
$GL(n,\C)$ of the form $g(\l)={\rm I}\ + \frac{P}{\l-z}$ for some $z\in \C$ and $P\in gl(n)$, we can
use residue calculus to choose a $gl(n)$-valued map $\ti P$ defined in an open neighborhood $\co$ of
$(0,0)$ in the $(\xi, \eta)$-plane so that 
 $$\ti E(\xi, \eta, \l)= g(\l) E(\xi, \eta,\l) \ti g(\xi, \eta, \l)^{-1}$$ is holomorphic in 
$\l\in \C\setminus \{0\}$ for each $(\xi, \eta)\in \co$, where $\ti g= I + \frac{\ti
P}{\l-z}$.  By a direct computation, one can see that 
 $$\ti E^{-1} \ti E_\xi= a\l + \ti u, \quad \ti E^{-1}\ti E_\eta= \l^{-1} \ti v$$
 for some $\ti u, \ti v$ defined on $\co$.  Hence $(a, \ti u, \ti v)$ is again a solution of the $-1$ equation.  We call 
 $$(a, u, v)\mapsto g\bu (a, u,v):= (a, \ti u, \ti v)$$
 a B\"acklund transformation of the $-1$ flow equation. 
We notice that the classical breather solutions of SGE can be
constructed from  the vacuum
solution $q=0$  by applying B\"acklund transformations twice with carefully placed poles.  Therefore
we can apply B\"acklund transformations
$2k$ times to construct $2k$-soliton breathers for SGE.  We show that corresponding wave maps
into $S^2$ are also homoclinic.   In fact, we generalize results of Shatah and Strauss to $1+1$ wave maps
into any compact symmetric spaces.

Note that if $a\in \cg$ is a constant, then $(a, 0,a)$ is a trivial solution of the $-1$ flow equation, and 
the $k$-soliton solutions of the $-1$ flow equation can be constructed by applying B\"acklund
transformations to it $k$ times.   If we choose $a$ so that $\exp(2\pi a)={\rm I}$ and place the poles of
the B\"acklund transformations carefully, then we can  obtain $k$-solitons of the
$-1$-flow equations, that are periodic in time or in space.  Such solutions are called {\it $k$-soliton
breathers\/}.  The wave maps into $G$ corresponding to $k$-soliton breathers are periodic in the space
variable. Wave maps from $S^1\times \R$ to $G$ are called {\it periodic wave maps\/} into $G$. 
The wave map corresponding to the trivial solution $(a,0,a)$ of the $-1$ flow equation is a stationary wave 
map into $G$, which is the geodesic $\g(x)=\exp(ax)$ in $G$.  We apply B\"acklund transformations to
these solutions to construct explicit $k$-soliton periodic wave maps.  Moreover, we compute the
asymptotic behavior of these periodic wave maps and prove that they are homoclinic.  We also
construct explicitly
$2k$-soliton homoclinic wave maps from
$S^1\times \R$ into $\C P^n$ and $4k$-soliton homoclinic wave maps into $S^{n-1}$.  

When $G$ is compact, a wave map from $\R^{1,1}$ to $G$ corresponding to a general $k$-soliton
solution of the
$-1$ flow equation usually oscillates as the space variable $|x|\to \infty$, i.e., does not have constant
boundary conditions at infinity.  To construct wave maps into $G$ that have good boundary conditions and
finite energy, we first note that  a wave map to a circle subgroup $T$ of $G$ is given essentially by a
solution of  the linear wave equation, so  wave maps to $T$ with finite energy and good boundary
conditions at $\pm
\infty$ can be written down easily, and B\"acklund transformations of such wave maps are again wave
maps having finite energy and constant boundary conditions at infinity.  

The Lax pair of the defocusing NLS 
$$q_t= \frac{i}{2}(q_{xx} - 2|q|^2 q)$$ satisfies the reality condition
coming from the non-compact Lie group $SU(1,1)$, and does not
have smooth solitons.  The theorem that B\"acklund transformations do not introduce singularities applies
only to the reality condition coming from a compact Lie group. Unfortunately there are many interesting
geometric problems in integrable systems for which solutions obtained via B\"acklund transformations do
have singularities.  Nevertheless, B\"acklund transformations can still be used in the non-compact case to
construct interesting examples as we will see next with
$G=SL(2,\R)$.

 It is known that the Cauchy problem for wave maps
from $\R^{1,1}$ to a complete Riemannian manifold
$N$ with smooth initial data in $L^2_1$ has long time existence (cf. \cite{Gu80}).  But this is no longer
true if we replace
$N$ by a pseudo-Riemannian manifold. There are counterexamples for $N=SL(2,\R)$ equipped with the
pseudo-Riemannian bi-invariant metric.

This paper is organized as follows.  In  section 2, we review   the
Lagrangian formulation of wave maps from $\R^{1,1}$ to $G$ and the corresponding Lax pair.   In section
3, we give the Hamiltonian formalism for wave maps  and compute the stable and unstable modes at
stationary solutions.  In section 4, we  prove the Lax pair of the $-1$ flow equation associated to $G$ is
gauge equivalent to the Lax pair of the equation for wave maps into $G$, and  give a bijection between
solutions of the $-1$ flow equation and wave maps
$s$ satisfying
$s(0,0)={\rm I\/}$.    In section 5, we review B\"acklund transformations for the $-1$ flow equation
associated to $SU(n)$, and apply these transformations to stationary wave maps to construct 
explicit $k$-soliton wave maps from $S^1\times \R$ to $SU(2)$.  In section 6, we prove the wave maps to
$SU(2)$ corresponding to $2k$-soliton breather solutions are homoclinic.  In section 7, we
explain the constraint condition for the $-1$ flow equation associated to $SU(2)$ so that the
corresponding wave maps into $SU(2)=S^3$ lie in 
$S^2$.  In section 8, we first recall a useful description of the compact symmetric space
$G/K$   imbedded as a totally geodesic submanifold in $G$, and then prescribe the
constraint condition for the $-1$ flow associated to $G$  that  insures that the corresponding wave
maps into
$G$ actually lie in a symmetric space. We emphasize the important case of wave maps into
$S^2=\frac{SU(2)}{SO(2)}$,  into $\C P^{n-1}$, and into
$S^{n-1}$.       Finally in section 9, we use B\"acklund transformations to construct examples of smooth Cauchy
data with constant boundary conditions at infinity and finite energy such that the Cauchy problem for wave
maps from
$\R^{1,1}$ to $SL(2,\R)$ have long time existence and also examples of initial data that develop singularities in
finite time.

\section{Wave map equation and its Lax pair}\label{fb}

We rewrite \eqref{am} as a first order system. 
Let $P= s^{-1}s_\xi$, and $Q= s^{-1}s_\eta$, i.e.,
\begin{equation} \label{bs}
s_\xi= s P, \quad s_\eta= sQ.
\end{equation}
The compatibility condition of the linear system \eqref{bs} gives 
\begin{align*}
&(s_\xi)_\eta = (sP)_\eta= s_\eta P + sP_\eta = sQP + sP_\eta= s(QP+ P_\eta)\cr
&\ \ = (s_\eta)_\xi = (sQ)_\xi = s_\xi Q + sQ_\xi = sPQ + sQ_\xi = s(PQ+ Q_\xi).
\end{align*}
This implies
$$QP+ P_\eta= PQ+ Q_\xi,$$
or equivalently, 
\begin{equation} \label{bq}
P_\eta- Q_\xi = PQ-QP= [P.Q].
\end{equation}
Combine equations \eqref{bq} and \eqref{am} to see that the wave map equation in characteristic
coordinates is 
$$P_\eta= -Q_\xi = \frac{1}{2}[P,Q],$$
i.e., 
\begin{equation} \label{bp}
(s^{-1}s_\xi)_\eta= -(s^{-1}s_\eta)_\xi = \frac{1}{2} \ [s^{-1}s_\xi, \ s^{-1} s_\eta].
\end{equation}
In other words, we have

\begin{prop} Let $(\xi, \eta)$ denote the light cone coordinate system of $\R^{1,1}$.  
 If $s:\R^{1.1}\to SU(n)$ is a wave map, then $A= \frac{1}{2}
s^{-1}s_\xi$ and
$B=\frac{1}{2}s^{-1}s_\eta$ satisfy the  first order system
\begin{equation} \label{br}
A_\eta= -B_\xi= [A,B].
\end{equation} 
Conversely, if $(A,B)$ is a solution of \eqref{br}, then there exists a unique
$s:\R^{1,1}\to SU(n)$ such that 
$$ s_\xi = 2sA, \quad s_\eta= 2sB, \quad s(0,0)=I.$$ 
Moreover,  $s$ satisfies \eqref{am}, i.e., $s$ is a wave map.   
\end{prop}

Next we formulate equation \eqref{br} as the condition for a family of connections to be
flat.   Recall that the curvature of a $gl(n)$-valued connection 
$$\left\{\frac{\p}{\p \xi}+P, \quad \frac{\p}{\p \eta}+Q\right\}$$
is defined to be 
$$F=\left[\frac{\p}{\p \xi} + P, \ \ \frac{\p}{\p \eta} +Q\right] = -P_\eta +Q_\xi + [P,Q].$$
The connection is {\it flat\/} if the curvature is zero.  
So the compatibility condition \eqref{bq}  for linear system \eqref{bs} is also the
condition for the connection $\left\{\frac{\p}{\p \xi}+P, \quad \frac{\p}{\p \eta}+Q\right\}$ to be
flat.  Another convenient way to write connection is as a  $gl(n)$-valued $1$-form 
$$\o= A\ d\xi + B\ d\eta.$$
Then the curvature is 
$$ d\o +\o \wedge \o = (-A_\eta+B_\xi +[A,B])\ d\xi\wedge d\eta.$$

It is easy to see that  the following statements are equivalent for  smooth maps $A, B:\R^2\to
gl(n)$:
\ben
\item  $\left[\frac{\p}{\p x} + A(x,t), \, \frac{\p}{ \p t} + 
B(x,t)\right] = 0$,
\item  The connection $1$-form $\o=Adx + Bdt$ is flat, i.e., 
$d\o=
-\o\wedge \o$. 
\item  $A_t-B_x = [A,B]$.
\item 
$ E_x= EA,\quad  E_t=EB,\quad E(0,0)=I$
 has a  unique solution $E:\R^2\to GL(n)$.    Such $E$ is called  {\it the trivialization\/} of
the flat connection $A\ dx + B\ dt$ (normalized at $(0,0)$).   
\een

  The wave map equation has a Lax pair (cf. \cite{ZakMik78, Uhl89}), i.e.,   there is  a one
parameter family of
$sl(n, \C)$-valued connection $1$-forms $\W_\l$
 on $\R^{1,1}$ defined in terms of $s:\R^{1,1}\to SU(n)$ and its
derivatives so that 
$\W_\l$ is flat for all $\l\in \C\setminus \{0\}$ if and only if $s$ satisfies \eqref{am}.   
We explain this next. 
Given $A, B:\R^{1,1}\to su(n)$ and $\l\in \C\setminus \{0\}$, consider the following
$gl(n,\C)$-valued connection $1$-form on $\R^{1,1}$:
$$\W_\l = (1-\l) A\ d\xi + (1-\l^{-1}) B\ d\eta.$$
We claim that $\W_\l$ is flat for  all $\l\in \C^*$ if and only if 
 $(A,B)$ is a solution of \eqref{br}.  
To see this, note that $\W_\l$ is flat is equivalent to
$$(1-\l)A_\eta- (1-\l^{-1})B_\xi = [(1-\l)A, (1-\l^{-1})B]= (2-\l-\l^{-1})[A,B]$$ 
for all $\l\in \C\setminus \{0\}$.  
Equate the coefficients of $\l^{-1}, \l$ and constant term to get
$$\bca
A_\eta=[A,B],&\cr A_\eta-B_\xi = 2[A,B],&\cr B_\xi= 
-[A,B].
\eca $$
 This is equivalent to \eqref{br},
 and we prove the claim.  We summarize our discussions:

\bprop \label{aa} (\cite{Uhl89, ZakMik78}).
 Let $s:\R^{1,1}\to SU(n)$ be a smooth map, $(\xi,\eta)$ the light cone coordinate system, and
 $$A=\frac{1}{2}(s^{-1}s_\xi), \quad B=\frac{1}{2}(s^{-1}s_\eta).$$  Then the following statements
are equivalent:
\ben
\item $s$ is a wave map.
\item $s$  is a solution of \eqref{am}.
\item $(A,B)$ is a solution of \eqref{br}.
\item The connection $1$-form 
 \begin{equation}\label{ab}
\W_\l= (1-\l)A\ d\xi + (1-\l^{-1}) B\ d\eta
\end{equation}
is flat  for all $\l\in \C\setminus\{0\}$,
\een 
\eprop

\bcor \label{ea}
If $A, B:\R^{1,1}\to su(n)$ satisfy equation \eqref{br}, then 
there exists $E(x,t,\l)$ such that 
 $$E^{-1}E_\xi=(1-\l)A, \quad E^{-1}E_\eta= (1-\l^{-1})B, \quad E(0,0,\l)=I$$
for all $\l\in \C\setminus\{0\}$, i.e., $E(\cdot, \cdot, \l)$ is the trivialization of the Lax pair $\W_\l$ defined by
\eqref{ab}.  Moreover, 
$s(\xi,\eta)= E(\xi, \eta,-1)$ is a wave map from $R^{1,1}$
to $SU(n)$,  $s^{-1}s_\xi= 2A$, and $s^{-1}s_\eta= 2B$.  
\ecor

A direct computation implies that the Lax pair $\W_\l$ of the wave map equation satisfies the following reality
condition:
\begin{equation}\label{cg}
\W_{\bar \l}^*+ \W_\l=0,
\end{equation}
where $\xi^*= \bar \xi^t$.  
We claim that the trivialization $E(x,t,\l)$ of $\W_\l$ satisfies the reality condition
\begin{equation} \label{ch}
E(x,t,\bar\l)^*E(x,t,\l)= I,
\end{equation} 
or equivalently,
$$E(x,t,\l)^{-1}= E(x,t,\bar \l)^*. $$
To see this, let $F(x,t,\l)=(E(x,t,\bar \l)^*)^{-1}$. Compute directly to get $F^{-1}dF= -\W_{\bar \l}^*$, which
is equal to $\W_\l$.  But $F(0,0,\l) ={\rm I\/}$.  Since both $E$ and $F$ satisfy the same linear differential
equation with the same initial condition, the uniqueness of ODE implies that $E=F$. This proves the claim. 

The $\l$ parameter seems redundant.  But it is this parameter that allows us to
 construct B\"acklund transformations and explicit solutions.  These will be explained in later sections.

\section{The Hamiltonian formulation of wave maps}\label{fc}

The functional $J$ defined by
\eqref{ev} for maps $s:S^1\times \R\to SU(n)$ is
$$J(s)=\frac{1}{2}\int_{\R^2}  ||s^{-1}s_t||^2 -||s^{-1}s_x||^2\ dx dt,$$
where $||y||^2= -\tr(y^2)$.  Viewed as a
functional on the space of curves from
$\R$ to
$C^{\infty}(S^1, SU(n))$, $J$ has two terms. The first term of $J$ is the kinetic energy and the
second term is the potential energy.  The Lagrangian formulation of the wave map equation  views the
equation as an equation for curves on the tangent bundle of
$\cm=C^\infty(S^1, SU(n))$.  In this section,  we use the Legendre
transformation to view the wave map equation as a Hamiltonian system on the cotangent bundle of
$\cm$, and compute the stable and unstable modes at stationary wave maps.  

Recall that  the cotangent bundle $T^*\cm$ of a manifold $\cm$ has a natural symplectic form
$w=d\tau$, where $\tau$ is the canonical 
$1$-form on $T^*\cm$ defined by 
$$\tau_\ell(v) = \ell(d\pi(v)),$$
where $\pi:T^*\cm\to \cm$ is the natural projection.  

Given a curve $\g:
(-\e, \e)\to C^\infty(S^1, SU(n))$ with $\g(0)=s$, we identify the tangent vector $\g'(0)$ as
$$(\g(0), \g(0)^{-1}\g'(0)) =(s, s^{-1}\d s).$$ 
 This identifies
$T\cm = \cm \times C^\infty(S^1, su(n))$.  Note that 
$$(v_1, v_2)= -\tr(v_1 v_2)$$
defines an inner product on $su(n)$.  So we can also identify $T^*\cm_s$ as
$T\cm_s$ via the
$L^2$ inner product:
$$\li v_1, v_2\ri = \int_0^{2\pi} -\tr(v_1v_2) dx.$$
   By definition of the canonical $1$-form on
$T^*\cm$, we get
$$\tau_{(s,v)} (s^{-1} \d s, \d v)= \li v, s^{-1} \d s\ri =\int_0^{2\pi} -\tr(v s^{-1}\d s)\ dx.$$
We use Cartan formula 
$$w(X,Y)= d\tau(X,Y) = X(\tau(Y)) - Y(\tau(X)) -\tau([X,Y])$$
to compute  the symplectic form $w=d\tau$ on $T^*\cm$.  Let $X(s, v)=(\eta_1, \d_1 v)$ and $Y(s,v)=
(\eta_2, \d_2 v)$ be two constant vector fields on $T^*\cm=T\cm$.  Then 
$$X(\tau (Y))= X\left(\int_0^{2\pi} -\tr(v(x)\eta_2(x))\ dx\right) = -\int_0^{2\pi} \tr(\d_1v, \eta_2) dx = \li
\d_1 v,
\eta_2\ri.$$  So we get
$$w_{(s,v)} ((\eta_1, \d_1 v), (\eta_2, \d_2 v))= \li \d_1 v, \eta_2\ri - \li \d_2 v, \eta_1\ri. $$
Consider the Hamiltonian $H:T^*\cm\to \R$, which is the sum of kinetic energy and potential energy,
i.e., 
$$H(s,v) = \frac{1}{2}\left(\li v, v\ri + \li s^{-1}s_x, s^{-1}s_x\ri\right) = -\frac{1}{2}\int_0^{2\pi}\tr(v^2 +
(s^{-1}s_x)^2) dx.$$
The Hamiltonian vector field $X_H$ of $H$ is the vector field satisfying 
$$dH_{(s,v)}(s^{-1}\d s, \d v) = w((s^{-1}\d
s,\d v, \ X_H(s,v))$$
for all $(s^{-1}\d s, \d v)$.  A direct computation shows that 
$$dH_{(s,v)}(s^{-1}\d s, \d v) = \li \d v, v\ri -\li s^{-1}\d s, (s^{-1}s_x)_x\ri.$$ 
So the Hamiltonian vector field for $H$ is 
$$X_{H}(s,v) = (v, (s^{-1}s_x)_x).$$ 
The Hamiltonian equation is 
$$s^{-1}s_t= v, \quad v_t=(s^{-1}s_x)_x,$$
which is the wave map equation \eqref{ae}.

\begin{prop} \label{au}
 The stationary points of  $X_H$ are $(s,0)$, where $s(x)=ce^{ax}$, $a\in su(2)$ a constant such that $e^{2\pi
a}=I$ and  $c\in SU(2)$ a constant.  
\end{prop}

\begin{proof}
 $X_H(s,v)=(0,0)$ if and only if $v=0$ and $(s^{-1}s_x)_x=0$. 
So $s^{-1}s_x=a$ for some constant $a\in su(n)$.  Hence $s(x)=ce^{ax}$
for some $c\in su(n)$.    Since
$s(2\pi)=s(0)$, $e^{2\pi a}=I$. 
\end{proof}

Note that stationary points of $X_H$ are closed geodesics of $SU(n)$.

Next we compute the linearization of $X_H$ at a stationary point, stable and unstable
subspaces.  We will do this calculation for $SU(2)$.  The calculations for other compact groups are
similar.    

Let $m$ be a non-zero integer,
$a=\diag(im, -im)$, and $s(x)=e^{ax}$.  The linearization of $X=X_H$ at the stationary point $(s,0)$
is
$$dX_{(s,0)}(s^{-1}\d s, \d v) = (\d v, \d(s^{-1}s_x)_x).$$ 
Set 
$$s^{-1}\d s= p, \quad \d v= q.$$
Compute directly to get
\begin{align*}
\d(s^{-1}s_x) &= -(s^{-1}\d s) s^{-1}s_x+ s^{-1}(\d s)_x \cr
&= -p a + s^{-1}(s p)_x = -p a + s^{-1}(s_x p+ sp_x)\cr
&= -p a + ap + p_x = p_x + [a, p].
\end{align*} 
So  
$$dX_{(s,0)} (p, q)= (q, p_{xx} + [a, p_x]).$$
The linearized equation is 
\begin{equation}\label{en}
\bca p_t= q, &\cr q_t= p_{xx}+[a,p_x].\eca
\end{equation}
The linearization of the wave map equation at $s$ is
\begin{equation} \label{co}
p_{tt}= p_{xx} + [a, p_x].
\end{equation}

We  compute the linear modes of the linear wave equation \eqref{en} next, i.e.,   solve the
following linear system for  $(p, q):S^1\to su(2)$:
$$\bca
q= kp,&\cr p_{xx}+[a,p_x]=kq.
\eca
$$
Substitute the first equation to the second to get
 \begin{equation}\label{bx}
p_{xx}+[a,p_x]-k^2p=0.
\end{equation}
Write \eqref{bx} in terms of entries of $p=(p_{ij})\in su(2)$ to get
\begin{equation} \label{bw}
\bca
(p_{11})_{xx} - k^2p_{11}=0, &\cr
(p_{12})_{xx} + 2im (p_{12})_x - k^2 p_{12}=0.
\eca
\end{equation}
This system is linear with constant coefficients. So it can be solved explicitly:
 \begin{align*}
p_{11}(x)&=\bca 
c_1+ c_2 x, & {\rm if}\  k=0,\cr
c_1 e^{kx} + c_2e^{-kx}, & {\rm if \ } k\not=0,
\eca\cr
p_{12}&=c_1 e^{\left(-im+\sqrt{k^2-m^2}\/\right)x} + c_2 e^{\left(-im-\sqrt{k^2-m^2}\/\right)x}
\end{align*}

We divide the computation into three cases:
\par

\ni {\bf (1)}   $k=0$.

Since in this case $p_{11}$ is linear and periodic with period
$2\pi$,  $p_{11}=c_1$ is a pure imaginary constant. 
Note $p_{12}= c_2 + c_3 e^{-2imx}$ with $c_2, c_3\in \C$. So the
nullity of $d(X_H)_{(s,0)}$ is $5$.  In fact, let $s_{cb}$ denote the stationary point $s_{c,b}(x)=
ce^{bx}$.  Then $\{s_{c,b}\n c\in SU(2), b\in su(2) \ {\rm is \ conjugate \ to \ } a\}$  is a five
dimensional stationary submanifold of $\cm$ and the tangent space at $s= s_{{\rm I\/}, a}$ is the kernel of
the linearization \eqref{co}.

\ni  {\bf (2)} $k\in \R\setminus\{0\}$.

 Note $p_{11}= c_1 e^{kx} + c_2e^{-kx}$ is periodic and $k$ is
non-zero and real implies that  $p_{11}=0$. 
 Since $p_{12}$ has period $2\pi$,  $m^2-k^2\geq
0$ and $\sqrt{m^2-k^2}$ is an integer.  So real non-zero eigenvalues of
$d(X_H)_{(s,0)}$ are
$$k=\pm \sqrt{m^2- j^2}, \quad 0\leq |j | < m, \ \ j \ {\rm integer\/}.$$
Eigenvectors for $k=\pm\sqrt{m^2-j^2}$ are $(p_k, q^\pm_k)$, where  
\begin{equation} \label{ci}
\begin{aligned}
p_k&= \bpm 0& c_1e^{-i(m+j)x} + c_2 e^{-i(m-j)x}\cr  -\bar c_1 e^{i(m+j)x} - \bar c_2
e^{i(m-j)x} &0\epm, \cr
 q^\pm_k&= \pm \sqrt{m^2-j^2}\  p_k.
\end{aligned}
\end{equation}

\ni {\bf (3)} $k\in \C\setminus\R$.

Since $p_{12}$ has to be periodic, 
$m^2-k^2>0$ and
$\sqrt{m^2-k^2}$ is an integer.  Hence $k^2$ must be real.  But $k$ is not real.  
So $k$ is pure imaginary.  In other  words,  $k= ic$ for some
$c\in \R$ and
$\sqrt{m^2+ c^2}$ is an integer.  Hence the non-real eigenvalues are
$$\pm i \sqrt{j^2-m^2}, \quad j > |m|, \ \  j \ {\rm integer\/}.$$

Recall that the {\it stable\/} ({\it unstable\/} resp.) subspace of $X_H$ at a stationary
point $(s,0)$ is the direct sum of the eigenspaces of $d(X_H)_{(s,0)}$ with eigenvalues
$k$ such that Re$(k)<0$ (Re$(k)>0$ resp.)  So the above computation gives

\begin{prop} \label{dj}
  Let $m$ be an integer, $a=\diag(im, -im)$, and $s(x)=e^{ax}$.  The
unstable subspace of the Hamiltonian vector field $X_H$ at the stationary point $(s,0)$ is 
$\oplus_{j=0}^{m-1} W^+_j$,
where $W^+_j$ is the eigenspace of $d(X_H)_{(s,0)}$ with eigenvalue $k=\sqrt{m^2-j^2}$
and is spanned by $(p_k, q_k^+)$ given in \eqref{ci}.
The stable subspace of the Hamiltonian vector field $X_H$ at $(s,0)$ is 
$\oplus_{j=0}^{m-1}W^-_j$,
where $W^-_j$ is the eigenspace of $d(X_H)_{(s,0)}$ with eigenvalue $k= -\sqrt{m^2-j^2}$
and is spanned by $(p_k, q_k^-)$ given in \eqref{ci}.
\end{prop}

\begin{cor}\label{cza}
Let $m$ be a positive integer, and $a=\diag(im, -im)$.  Then the linearization of the wave map equation at
the stationary wave map $s(x,t)= e^{ax}$ is 
$$\xi_{tt}= \xi_{xx} + [a, \xi_x].$$
Moreover, the stable and unstable modes corresponding to $\pm\sqrt{m^2-j^2}$ are 
\begin{subequations}\label{da}
\begin{gather}
p_{m,j}^-(x,t)= -e^{-\sqrt{m^2-j^2}\ t} \bpm 0& ce^{i(-m\pm j)x}\cr -\bar c e^{-i(-m\pm j)x} &0
\epm \label{da1}\\
p_{m, j}^+(x,t)= e^{\sqrt{m^2-j^2}\ t} \bpm 0& ce^{i(-m\pm j)x}\cr -\bar c e^{-i(-m\pm j)x} &0
\epm \label{da2}
\end{gather}
\end{subequations}
respectively, where $c\in \C$ is a constant, $j$ is an integer and $|j|<m$.  
\end{cor}

\section{ The $-1$ flow equation and the wave map}\label{fd}

We  give  a correspondence between solutions of the $-1$ flow equation  \eqref{eu}
and wave maps. 

In characteristic coordinate $(\xi, \eta)$, the $-1$ flow equation \eqref{eu} associated to $SU(n)$ is the
following  system for $(a,u,v):\R^2\to \Pi_{i=1}^3 su(n)$:
\begin{equation} \label{ac}
\bca 
a_\eta= 0, &\cr
u_\eta = [a, v],&\cr v_\xi = -[u,v]. 
\eca
\end{equation}

A direct computation implies that 

\begin{prop} \label{ag}
 The map $(a,u, v):\R^2\to \Pi_{i=1}^3 su(n)$ is a solution 
of the $-1$-flow equation
\eqref{ac}  associated to $SU(n)$ if and only if 
\begin{equation} \label{ap}
\o_\l= (a\l +u)\ d\xi + \l^{-1} v\ d\eta
\end{equation}
is flat for all $\l\in
\C\setminus \{0\}$. 
\end{prop}

Note that the Lax pair $\o_\l$ of the above $-1$ flow equation satisfies the reality condition
\eqref{cg}. So the trivialization $E(x,t,\l)$ of $\o_\l$ satisfies the reality condition \eqref{ch}.

Recall that the gauge transformation of $g:\R^{1,1}\to GL(n,\C)$ of the 
connection
\begin{equation*}
\left\{\frac{\p}{ \p \xi} + A, \quad \frac{\p}{ \p \eta}+ B\right\}
\end{equation*}
is $$\left\{g\left(\frac{\p}{ \p \xi} + A\right)g^{-1}, \quad
g\left(\frac{\p}{\p \eta}+ B\right)g^{-1}\right\}.$$ Direct computation gives
\begin{align*}
g\left(\frac{\p}{ \p \xi} + A\right)g^{-1}&=\frac{\p}{ \p \xi} +
gAg^{-1} -g_\xi g^{-1} , \cr
g\left(\frac{\p}{\p \eta}+ B\right)g^{-1} &=\frac{\p}{ \p \eta} +
gBg^{-1}-g_\eta g^{-1}.
\end{align*}
Since 
$$\left[g\left(\frac{\p}{ \p \xi} + A\right)g^{-1}, \ \  g\left(\frac{\p}{\p \eta}+ B\right)g^{-1} \right]
= g\left[\frac{\p}{ \p \xi} + A, \ \ \frac{\p}{\p \eta}+ B \right] g^{-1},$$ the gauge transformation of a
flat connection is again flat.   Written in terms of connection $1$-form 
$\o= A \ d\xi + B\ d\eta$, the gauge transformation $g\ast \o$ is 
$$g\ast \o = g\o g^{-1} - dg g^{-1}.$$

It is easy to check that if $E$ is the trivialization of the
flat connection $\o$, 
 then $g(0,0)^{-1}Eg^{-1}$ is the trivialization of $g\ast \o$.

Below we show that the Lax pairs of wave map equation and the $-1$ flow equation are gauge
equivalent and give a correspondence between wave maps and solutions of the $-1$ flow equation. 

\begin{thm}\label{ah}
  Let $(a,u,v)$ be a solution to the $-1$-flow
equation \eqref{ac} associated to $SU(n)$, and $\Phi(\xi,\eta, \l)$ the trivialization 
of
$$\o_\l = (a(\xi,\eta)\l + u(\xi,\eta))\ d\xi + \l^{-1}v(\xi,\eta)\ d\eta.$$  
Set
$\Phi(\l)(\xi,\eta) =\Phi (\xi,\eta,\l)$. Then
$s=\Phi(-1)\Phi(1)^{-1}$ is a wave map from $\R^{1,1}$ to $SU(n)$,  and
\begin{equation} \label{ec}
s^{-1}s_\xi=-2\Phi(1) a \Phi(1)^{-1}, \quad s^{-1} s_\eta= -2 \Phi(1) v\Phi(1)^{-1}.
\end{equation}
  Conversely,
suppose $s:\R^{1,1}\to SU(n)$ is a wave map and $s(0,0)={\rm I}$.  Let   $\psi(\xi,\eta)$ be  the
solution of $\psi^{-1}\psi_\eta = \frac{1}{ 2} s^{-1}s_\eta$ with 
$\psi(\xi,0)={\rm I}$, and  
\begin{equation} \label{bn}
\bca
a(\xi, \eta)= \frac{1}{ 2}(s^{-1}s_\xi)(\xi, 0), &\cr
u(\xi,\eta)= a(\xi) -(\psi_\xi\psi^{-1})(\xi,\eta), &\cr 
v= -\frac{1}{ 2}
\psi s^{-1}s_\eta\psi^{-1}.
\eca
 \end{equation}
Then  $(-a,u,v)$ is a solution of the 
$-1$ flow equation associated to $SU(n)$ and $s$ is the wave map corresponding to $(-a, u,v)$.  
\end{thm}

\begin{proof}
 A direct computation gives $$\Theta_\l= \Phi(1)\ast \o_\l = (1-\l)(-\Phi(1)
a\Phi(1)^{-1}) d\xi + (1-\l^{-1})(-\Phi(1) v\Phi(1)^{-1})d\eta.$$ The 
trivialization of
$\Theta_\l$ is $\Phi(\l)\Phi(1)^{-1}$.  By 
Corollary \ref{aa}, $\Phi(-1)\Phi(1)^{-1}$ is a wave map. 

To prove the converse, 
 set $A=\frac{1}{ 2} s^{-1}s_\xi$ and $B=\frac{1}{ 2} s^{-1}s_\eta$. 
By Proposition \ref{aa}, we have
$A_\eta= -B_\xi = [A,B]$.
Set $\tilde A(\xi,\eta)=\psi(\xi,\eta)^{-1}a(\xi)\psi(\xi,\eta)$. 
A direct
computation implies that $\tilde A_\eta= [\tilde A, B]$.  But $A$ 
satisfies the same
differential equation as $\tilde A$, i.e., $A_\eta= [A, B]$, and  
$$\tilde A(\xi,0)=A(\xi, 0)=a(\xi) = \frac{1}{ 2} (s^{-1}s_\xi)(\xi, 
0).$$ By the
uniqueness of solutions of ordinary differential equation we have
$$A(\xi,\eta)=\tilde A(\xi, \eta)= \psi^{-1}(\xi,\eta) a(\xi) 
\psi(\xi,\eta).$$  

Apply gauge transformation of $\psi$ to the Lax pair of the wave map 
$$\W_\l = (1-\l) A \ d\xi + (1-\l^{-1}) B\ d\eta$$ to get 
\begin{equation*}
\psi\ast \W_\l = (-a\l + a-\psi_\xi \psi^{-1}) \ d\xi -\l^{-1}\psi B\psi^{-1}\ d\eta.
\end{equation*}
Since $\W_\l$ is flat. so is $\psi\ast \W_\l$. 
  It follows from Proposition \ref{ag} that $(-a, u,v)$ is a solution of  the $-1$ flow,  where $a, u, v$ are
defined by \eqref{bn}. 

Let $F(\xi,\eta,\l)$ denote the trivialization of $\W_\l$. Since $\W_1=0$, $F(\xi,\eta,
1)$ is a constant.  But $F(0, 0, \l)= {\rm I\/}$.  Thus $F(\xi,\eta, 1)={\rm I\/}$.  It
follows from  Corollary \ref{ea} that the harmonic map $s(\xi,\eta)=F(\xi, \eta,-1)$.  
The trivialization of $\psi\ast\W_\l$ is 
$$E(\xi,\eta,\l)= F(\xi,\eta,\l) \psi(\xi,\eta).$$  But we have proved $\psi\ast \W_\l$ is
the Lax pair for $(-a, u,v)$.  So the wave map corresponding to $(-a, u,v)$ is 
$$E(\xi,\eta, -1)E(\xi,\eta,1)^{-1} = F(\xi,\eta, -1)F(\xi,\eta, 1)^{-1} = F(\xi, \eta,
-1)=s.$$
\end{proof}

The proof of the above Theorem implies that the Lax pair \eqref{ab} of the wave map equation is
gauge equivalent to the Lax pair \eqref{ap} of the $-1$ flow equation. 

\beg \label{at}  {\bf closed geodesics}
\par 
 
 Let $a=\diag(im_1, \cdots, im_n)\in su(n)$, where
$2m_1, \cdots,
2m_n$ are integers.  Then  $(a,0,a)$ is a 
 solution of the $-1$
flow equation, the corresponding Lax pair is
$$\o_\l= a\l \ d\xi + a\l^{-1} d\eta,$$
and the trivialization of $\o_\l$ is
$$E_0(\xi,\eta, \l)= e^{a(\l\xi + \l^{-1}\eta)}.$$
The corresponding wave map constructed in Theorem \ref{ah} is the stationary wave map
$$s_0(\xi,\eta)=E_0(\xi,\eta, -1)E_0(\xi,\eta, 1)^{-1} =
e^{-2ax}= \diag(e^{-2im_1 x},\ldots, e^{-2im_n x}),$$
a closed geodesic in $SU(n)$. 
\eeg

Note that $SU(2)$ equipped with the bi-invariant metric is isometric to
the standard
$S^3$ because
$$SU(2)=\left\{\bpm z& -\bar w\cr w& \bar z\epm \ \bigg| \   z,w\in \C, \ 
|z|^2+|w|^2=1\right\}$$
is isometric to $S^3$ in $\C^2= \R^4$ via 
$$\bpm z& -\bar w\cr w& \bar z\epm\ \mapsto\  \bpm z\cr w\epm.$$

\beg \label{ed} {\bf wave maps into a great circle\/}
\par
Let $h(\xi)$ and $k(\eta)$ be smooth real valued functions on $\R$,
$u=0$, and 
$$a=h'(\xi)\bpm i&0\cr 0&-i\epm, \quad b= k'(\eta)\bpm i&0\cr 0&-i\epm.$$
Then $(a, 0, b)$ is a solution of the $-1$ flow equation associated to $SU(2)$, its Lax pair
is $\o_\l= a(\xi) \l\ d\xi + b(\eta)\l^{-1} d\eta$, and its trivialization is
$$E(\xi, \eta, \l)= \bpm e^{i(h(\xi)\l + k(\eta)\l^{-1})} &0\cr 0 & e^{-i(h(\xi)\l +
k(\eta)\l^{-1})}\epm.$$
The wave map corresponding to $(a, 0, b)$ is 
$$s(\xi,\eta)= \diag(e^{-2i(h(\xi) + k(\eta))},\  e^{2i(h(\xi) + k(\eta))}),$$
which lies in the great circle $|z|=1$ and $w=0$ in 
$$SU(2)= \left\{ \bpm z&-\bar w\cr w& \bar z\epm\bigg| z, w\in \C, \
|z|^2+|w|^2=1\right\} = S^3.$$
Hence it is also a wave map into the circle $S^1$.   The equation for wave maps into
$S^1$ is essentially the linear wave equation, and general solutions are of the form
$h(\xi)+ k(\eta)$.  Note that  if $h$ and $k$ decay at $\pm \infty$, then the
corresponding wave map tends to ${\rm I\/}$ as $|x|\to \infty$.  
\eeg

\section{B\"acklund transformations}\label{fe}

In this section, we use B\"acklund transformations to construct $k$-soliton solutions of the $-1$ flow equation, and use Theorem
\ref{ah} to construct the corresponding $k$-soliton wave maps.  Most of these wave maps oscillates as the space
variable $x$ tends to $\pm \infty$, but some of these wave maps are periodic in $x$.  Note that wave  maps into a
great circle of $SU(2)$ can be written in terms of solutions of the linear wave equation.  We show that if $s$ is a wave
map into a great circle so that $s$ has constant boundary condition at $\pm\infty$ and finite energy, then the
new wave maps obtained by applying B\"acklund transformations to $s$ also have constant boundary condition
and finite energy.

First we review the construction of B\"acklund transformations of the $-1$ flow equation.  Let $(a,u,v)$ be a
solution of the $-1$ flow equation \eqref{ac} associated to $SU(n)$, and $E(x,t,\l)$ the trivialization of its Lax pair \eqref{ap}
$\o_\l$, i.e., 
\begin{equation}\label{cz}
\bca E^{-1}E_\xi= a\l + u, \cr
E^{-1}E_\eta= \l^{-1} v,\cr
E(0,0,\l)= {\rm I\/}.
\eca
\end{equation}
Since the right hand side of \eqref{cz} is holomorphic in parameter $\l\in \C\setminus 0$, the solution
$E(x,t,\l)$ is holomoprhic in
$\l\in \C\setminus 0$.  Because $\o_\l$ satisfies the reality condition \eqref{cg}, $E$ satisfies \eqref{ch}.  

    Let $\pi$ be a Hermitian projection of $\C^n$ onto a
complex linear  subspace $V$, $\pi^\perp=I-\pi$ the 
projection onto the orthogonal
complement $V^\perp$,  and $z\in \C$.
Let $g_{z,\pi}:\C\to GL(n, \C)$ denote the rational map defined by 
\begin{equation}\label{do}
g_{z,\pi}(\l)=\pi+ \frac{\l-z}{ \l-\bar z} 
\pi^\perp = I + \frac{\bar z-z}{
\l-\bar z} \pi^\perp.
\end{equation}  
We call $g_{z,\pi}$ a {\it simple element\/}.   
A direct computation shows that $g_{z,\pi}$ satisfies the reality condition \eqref{ch}:
$$(g_{z,\pi}(\bar\l))^*g_{z,\pi}(\l)=I.$$  
In particular,
	$$g_{z,\pi}^{-1}(\l)= (g_{ z, \pi}(\l))^*= \pi + \frac{\l-\bar z}{\l-z} \ \pi^\perp = g_{\bar z, \pi}(\l).$$

To construct B\"acklund transformations for the $-1$ flow equation, 
we first find $\ti E$ and $\ti g$ so that 
$$g_{z,\pi}(\l) E(x,t,\l) = \ti E(x,t,\l) \ti g(x,t,\l)$$
with $\ti E$ holomorphic in $\l\in \C\setminus 0$ and $\ti g$  holomorphic in a neighborhood of
$\{ 0, \infty\}$.  Since the left hand side has a pole at $\l=\bar z$, so $\ti g$ must have too. 
In fact, $\ti g$ can be taken to be the form $g_{z, \ti\pi(x,t)}(\l)$ for some projection
$\ti\pi(x,t)$.  Moreover,  $\ti E$ is the trivialization of a new solution of the $-1$ flow equation. 
We state the results more precisely below.

\begin{thm} \label{ak} (\cite{TerUhl00a}).  
Let $(a,u,v)$ be a solution of the $-1$ flow
equation \eqref{ac}, and $E(\xi,\eta,\l)$ the trivialization of the corresponding Lax pair $\o_\l$, i.e., 
$$E^{-1}E_\xi= a\l + u, \quad E^{-1} E_\eta= \l^{-1} v, \quad E(0,0,\l) = I.$$ 
Let $z\in  \C\setminus \R$, and $\pi$ the
projection onto a linear subspace $V$ of $\C^n$.   Set 
\begin{align*}
&\tilde V(\xi,\eta) = E(\xi,\eta,z)^*(V) \cr 
 &\tilde \pi(\xi,\eta)=\, {\rm the\, Hermitian\ projection\, of\  }\C^n\ {\rm onto\, } \tilde
V(\xi,\eta)\cr 
&\tilde u = u + (z-\bar z)[\tilde \pi, a],\cr
&\tilde v = \frac{1}{ \n z\n^2}(\bar z\tilde \pi + z\tilde
\pi^\perp) v (z\tilde \pi + \bar z\tilde \pi^\perp),\cr 
&\tilde  E(\xi,\eta,\l) = g_{z,\pi}(\l) E(\xi,\eta, \l) g_{z,\tilde\pi}(\l)^{-1}
\cr
&\qquad = \left(\pi + \frac{\l-z}{
\l-\bar z} \pi^\perp\right) E(\xi,\eta,\l) \left(\tilde \pi(\xi,\eta) + 
\frac{\l-\bar z}{ \l-z} \tilde
\pi^\perp(\xi,\eta)\right).
\end{align*} 
Then 
\ben
\item $\ti E(\xi, \eta, \l)$ is holomorphic for $\l\in \C\setminus 0$,
\item $(a, \tilde u, \tilde v)$ is a new solution of the $-1$ flow equation,
\item  $\ti E$ satisfies the reality condition \eqref{ch}, and $\tilde E$ is the trivialization of the Lax pair of $(a, \tilde u, \ti v)$. 
\een
\end{thm}

We sketch the proof of this Theorem.  Let $g= g_{z,\pi}$, and $\ti g= g_{z,\ti\pi(x,t)}$.  Note that 
$\ti E=gE\ti g^{-1}$ is holomoprhic for $\l\in \C\setminus \{0, z, \bar z\}$ and has 
poles at $\l= z$ and
$\bar z$ of order $\leq 1$.  Use definition of $\ti \pi$ to prove that the residues of $\ti E$ at $\l=z$ and
$\l=\bar z$ are zero.  Thus $\ti E$ is holomorphic for $\l\in \C\setminus 0$.  Let $\ti \o_\l= \ti E^{-1} d\ti E$.  Then
\begin{equation}\label{dx}
\ti \o_\l = \ti g \o_\l \ti g^{-1}- (d\ti g) \ti g^{-1} = \ti g\ast \o_\l.
\end{equation}
Expand $\ti E^{-1}\ti E_\xi$ in $\l$ to see that its leading term is $a\l$.  Since $\ti E^{-1}\ti E_\xi$ is
holomorphic in $\l\in \C$, it must be of the form $a\l + \ti u$ for some $\ti u$.  A similar argument
implies that $\ti E^{-1} \ti E_\eta$ must be of the form $\l^{-1} \ti v$.  But $\ti \o_\l= (a\l + \ti u)\ d\xi
+ \l^{-1} \ti v \ d\eta$ is flat.  So $(a,\ti u, \ti v)$ is a solution of the $-1$ flow equation.  The
formula of $\ti u, \ti v$ can be computed from \eqref{dx}.

Theorem \ref{ak} gives an algebraic method to construct new solutions from a given solution of the $-1$
flow equation if the trivialization of the Lax pair of the given solution is known.  Let $g_{z,\pi}\bu (a,u,v)$ denote the solution
$(a, \tilde u, \tilde v)$ constructed in Theorem \ref{ak}, and the transformation $(a, u,v)\mapsto
g_{z,\pi}\bu (a,u,v)$ is called a {\it B\"acklund transformation\/} of the $-1$ flow equation.   Let
$s$ be the wave map corresponding to
$(a,u,v)$ given by Theorem \ref{ah}, and $g_{z,\pi}\bu s$ the wave map corresponding to
$g_{z,\pi}\bu (a, u, v)$.  We call $s\mapsto g_{z,\pi}\bu s$ a B\"acklund
transformation of wave maps.

In the next two examples, we use B\"acklund transformations to construct explicit wave maps
into $SU(2)$.  

\begin{eg} \label{bt} {\bf periodic $1$-soliton wave map\/}
\hfil

 Let $2m>0$ be an integer, $a=\diag(im,-im)$. 
We have seen in Example \ref{at} that $(a,0,a)$ is a solution of the $-1$ flow equation, its Lax pair is
$$\o_\l= a\l \ d\xi + a\l^{-1} \ d\eta,$$
its trivialization is $E_0(\xi, \eta, \l) = e^{a\l\xi + a\l^{-1} \eta}$, 
and the corresponding wave map is the stationary wave map $s_0(x, t)= e^{-2ax}$.   
Since we have the trivialization for $(a,0,a)$, we can apply B\"acklund transformation to $(a,0,a)$.  Let 
$z=e^{i\o}$,
$q_0=\bpm 1\cr 1\epm$, and $\pi$ the Hermitian projection onto $V=\C q_0$.   We use  Theorem \ref{ak} to
compute $g_{z,\pi}\bu (a,0,a)$ next.   First we get
\begin{align*}
\tilde q(\xi,\eta) &= \exp({a(z\xi+ z^{-1}\eta)})^*q_0= \exp(-a(\bar z \xi + \bar
z^{-1}\eta))q_0\cr &= \exp(-a(e^{-i\o} \xi + e^{i\o} \eta))q_0\cr 
& = \exp(-a(\cos \o (\xi+\eta) +
i\sin\o (-\xi+
\eta)))q_0\cr &= \exp(-a (x\cos\o - it\sin\o)\bpm 1\cr 1\epm
=\bpm e^{-(i mx\cos\o + mt\sin\o)}\cr  e^{imx\cos\o + mt\sin\o} \epm.
\end{align*}
Therefore the projection $\ti \pi(x,t)$ of $\C^2$ onto $\C \ti q(x,t)$ is
$$\ti \pi(x,t)= \frac{1}{e^{2mt\sin \o} + e^{-2mt\sin\o}}
\bpm  e^{-2mt\sin \o} & e^{-2i mx\cos\o} \cr e^{2imx\cos\o} & e^{2mt\sin\o}
\epm,
$$
the trivialization of $g_{z,\pi}\bu (a,0,a)$ is
\begin{align*}
E_1(\xi,\eta, \l)&= g_{z,\pi}(\l) E_0(\xi,\eta,\l) g_{z.\tilde\pi(\xi,\eta)}(\l)^{-1}\cr
&= \left(\pi + \frac{\l-z}{\l-\bar z} \pi^\perp\right) e^{a(\l\xi + \l^{-1}\eta)} \left(\ti \pi(x,t)
+\frac{\l-\bar z}{\l-z}
\ti\pi^\perp(x,t)\right).
\end{align*}
By Theorem \ref{ah},  the wave map $s$ corresponding to 
$g_{z,\pi}\bu (a,0,a)$ is 
\begin{align*}
& s(x,t) = E_1(x,t,-1) E_1(x,t, 1)^{-1}\cr
&= g_{z,\pi}(-1) E_0(x,t,-1) g_{z,\ti\pi(x,t)}(-1)^{-1} g_{z,\ti\pi(x,t)}(1) E_0(x,t,1)^{-1}
g_{z,\pi}(1)^{-1}.
\end{align*}
A direct computation gives
\begin{align*}
&E_0(x,t,-1)= E_0(x,t,1)^{-1} = e^{-ax},\cr
&g_{z,\ti\pi(x,t)}(-1)^{-1} g_{z,\ti\pi(x,t)}(1)= \ti\pi(x,t)-\ti\pi(x,t)^\perp.
\end{align*}
Hence
\begin{equation}\label{cj}
s(x,t)=g_{z,\pi}(-1) e^{-ax} (\ti \pi (x,t)- \ti \pi (x,t)^\perp) e^{-ax} g_{z,\pi}(1)^{-1}.
\end{equation}
In particular,  the first column of $e^{-ax} (\ti \pi(x,t)-\ti\pi(x,t)^\perp)e^{-ax}$
$$S(x,t)=\bpm - e^{-2imx}\tanh (2mt\sin \o)\cr e^{2i mx \cos\o} 
\sech(2mt\sin \o)\epm$$ is a
wave map into $S^3$. 

Note that if $\cos\o= \frac{j}{2m}$ for some integer $j$, then  $s$ is periodic in
$x$ with period $2\pi$.  In this case,  $S$ is a wave map from $S^1\times \R$ to $S^3$. 

\end{eg}

\beg  {\bf $1$-soliton wave map into $SU(2)$\/}
\hfill

  Let $z=r+is$, and $a$, $\pi$ as in Example \ref{bt}.  We derive the formula for $\ti s= g_{z,\pi}\bu s_0$.  
A direct computation as in Example \ref{bt} implies that 
$$\ti\pi(\xi, \eta)= \frac{1}{1+e^{2A}}\bpm 1& e^{A-iB}\cr e^{A+iB}& e^{2A}\epm,$$
where 
\begin{align*}
A&= ms\left(x+t -\frac{x-t}{|z|^2}\right) = ms\left(\left(1-\frac{1}{|z|^2}\right) x + \left(1+\frac{1}{|z|^2}\right)
t\right),\cr
 B&=m r\left(x+t+\frac{x-t}{|z|^2}\right) = mr\left(\left(1+\frac{1}{|z|^2}\right)x +
\left(1-\frac{1}{|z|^2}\right)t
\right).
\end{align*}
The first column  of $E_0(x,t,-1)(2\ti\pi(x,t)-{\rm I\/}) E_0(x,t,1)^{-1}$ is
\begin{equation*}
\ti S = \bpm -e^{-2imx}\tanh A \cr e^{iB} \sech A\epm,
\end{equation*}
which is a wave map into $S^3$. 
Note that $\ti S$ is periodic in $x$ if $|z|=1$ and $r$ is rational, and oscillates as $|x|\to \infty$ if $|z|\not=1$.  
\eeg

\beg \label{bk} {\bf $k$-soliton wave map from $S^1\times \R$ to $SU(2)$\/}
\hfil

We apply B\"acklund transformations $k$ times to construct $k$-soliton wave maps from
$S^1\times
\R$ to $SU(2)$.  Let $a, q_0, V$, and $\pi$ be as in Example \ref{bt}.  Let
$$z_j=e^{i\o_j}= \frac{r_j+i \mu_j}{ m}, \qquad 2r_j \ {\rm an \ 
integer,\ }\  | r_j| < m, \quad \mu_j= \sqrt{m^2-r_j^2}$$ for $j=1, \cdots, k$.  Set 
$$(a,u_j,v_j) = g_{z_j,\pi}\bu(g_{z_{j-1},\pi}\bu \cdots \bu 
g_{z_1,\pi}\bu  (a,0,a)\cdots),$$ 
$E_j$ the trivialization of the Lax pair of $(a, u_j,v_j)$, 
$$q_j(x,t)= E_{j-1}(x,t,z_j)^*(q_0)= E_{j-1}(x,t,\bar z_j)^{-1}(q_0),$$
$\pi_j(x,t)$ the Hermitian projection of $\C^n$ onto $\C q_j(x,t)$, i.e., 
$$\pi_j(x,t)=\frac{ q_j(x,t)q_j^*(x,t)}{||q_j(x,t)||^2}.$$
Set 
$$g_j(\l) = g_{z_j, \pi}(\l), \quad \ti g_j(x,t,\l)= g_{z_j, \pi_j(x,t)}(\l).$$
By Theorem \ref{ak} and induction, we have 
\begin{subequations}\label{cf}
\begin{gather} 
E_j(x,t,\l)= (g_j\cdots g_1)(\l) E_0(x,t,\l) (\ti g_j\cdots \ti g_1)^{-1}(x,t,\l), \label{cf1}\\
q_j(x,t)= (\ti g_{j-1}\cdots \ti g_1)(\bar z_j) E_0(x,t,z_j)^*(g_{j-1}\cdots g_1)^*(z_j)(q_0),\label{cf2}
\end{gather}
\end{subequations}
where $E_0(x,t,\l)= e^{a(\l \xi + \l^{-1}\eta)}$.

Since $\pi(q_0)=q_0$, $\pi^\perp(q_0)=0$,  
\begin{equation} \label{es}
g_i(z_j)^*(q_0 )= \left(\pi+ \frac{\bar z_j-\bar z_i}{\bar z_j- z_i}\pi^\perp\right)(q_0)= q_0.
\end{equation}
Hence $(g_j\cdots g_1)^*(z_j)(q_0) =  q_0$, and
\begin{align}\label{bz}
q_j(x,t) &= (\ti g_{j-1}\cdots \ti g_1)(\bar z_j)E_0(x,t,z_j)^* (q_0)\cr
&= (\ti g_{j-1}\cdots \ti g_1)(\bar z_j)\bpm e^{-(ir_j x+ \mu_j t)} \cr e^{ir_jx+ \mu_j
t}\epm.
\end{align} 

The  wave map 
corresponding to $(a, u_j, v_j)$ constructed in Theorem \ref{ah} is
$$s_j(x,t)=E_j(x,t,-1) E_j(x,t,1)^{-1}.$$

From Example \ref{bt}, we get
\begin{subequations} \label{by}
\begin{gather}
q_1(x,t)= \bpm e^{-(ir_1x + \mu_1 t)} \cr e^{ir_1 x+ \mu_1 t}\epm, \label{by1} \\
\pi_1(x,t)= \frac{1}{e^{2\mu_1 t}+ e^{-2\mu_1t}} \bpm 
e^{-2\mu_1 t}& e^{-2ir_1 x}\cr 
e^{2ir_1 x}& e^{2\mu_1 t}
\epm,\label{by2} \\
E_1(x,t,\l)= g_1(\l) e^{a(\l\xi+ \l^{-1}\eta)} \ti g_1(x,t,\l)^{-1}, \label{by3}\\
s_1(x,t) = g_1(-1) e^{-ax} (\pi_1(x,t)-\pi_1^\perp(x,t)) e^{-ax} g_1(1)^{-1}.\label{by4}
\end{gather}
\end{subequations}
  Since $r_1$ is an integer, $\pi_1, \ti g_1,  s_1$ are periodic in $x$ with
period $2\pi$.    

Note that 
\begin{equation}\label{et}
E_0(x,t,z_j)^*(q_0)= \bpm e^{-(ir_j x+ \mu_j t)}\cr e^{ir_jx+ \mu_j t}\epm.
\end{equation}
Use \eqref{cf2}, \eqref{es} and \eqref{et} to get
$$q_2(x,t)= \ti g_1(x,t, \bar z_2) \bpm e^{-(ir_2 x + \mu_2 t)}\cr e^{ir_2 x + \mu_2 t}\epm.$$
Since $\ti g_1$ is periodic in $x$ and $r_2$ is an integer, $q_2$ is periodic in $x$ with period $2\pi$.  So is $\pi_2$.  
The wave map corresponding to $(a,u_2, v_2)$ is
\begin{equation*}
s_2(x,t) = b_2e^{-ax} \ti g_1(x,t,-1)^{-1}
(\pi_2(x,t)-\pi_2(x,t)^\perp) \ti g_1(x,t,1) e^{-ax}c_2^{-1},
\end{equation*}
where $b_2= (g_2g_1)(-1)$, $c_2=(g_2g_1)(1)$, and $\pi_2(x,t)$ is the projection onto 
$$q_2(x,t)= \left(\pi_1(x,t) + \frac{\bar z_1-z_1}{\bar z_2 -\bar z_1}\pi_1^\perp(x,t)\right)\bpm 1\cr
e^{2ir_2x +\mu_2 t}\epm.$$
  By assumption $2r_2$ is an integer.  So $q_2$ and $s_2$ are periodic in $x$ with period $2\pi$.

It now follows from induction that
\ben 
\item the wave map corresponding to $(a, u_j, v_j)$ is
\begin{equation}\label{ca}
s_j(x,t)= b_j e^{-ax} (\ti g_1^{-1}\cdots \ti g_j^{-1})(x,t,-1) (\ti g_j\cdots \ti g_1)(x,t,1) c_j^{-1},
\end{equation} 
where 
\begin{equation} \label{ck}
b_j= (g_j\cdots g_1)(-1), \quad c_j=(g_j\cdots g_1)(1),
\end{equation}
\item $s_j$ is periodic in $x$ with period $2\pi$, i.e., $s_j$ is a wave map from $S^1\times \R$ to
$SU(2)$.   
\een
\eeg

Next we compute the wave map constructed by applying B\"acklund transformation to wave maps into $S^1$
constructed in Example \ref{ed}.  

\beg \label{ee} 
Let $z=r+is$, $v_0=\bpm 1\cr 1\epm$, $h, h$, $a, b$ and $s$ as in Example \ref{ed}. 
We use Theorem \ref{ak} to construct $(a, u, v)=g_{z,\pi}\ast (a,0,b)$.   A direct computation gives
\begin{equation*}
q(\xi,\eta)= E(\xi, \eta,z)^*(v_0)= \bpm e^{-i(h(\xi)\bar z + k(\eta)\bar z^{-1})}\cr 
e^{i(h(\xi)\bar z + k(\eta)\bar z^{-1})}\epm,
\end{equation*}
which is parallel to 
$$\bpm 1\cr e^{2s(h(\xi)-\frac{k(\eta)}{r^2+s^2})} e^{2i r(h(\xi)+\frac{k(\eta)}{r^2+s^2})}\epm.$$
So the projection $\ti\pi(\xi, \eta)$ onto $\C q(\xi,\eta)$ is
$$\ti\pi(\xi,\eta)= \frac{1}{1+ e^{4s(h(\xi)-\frac{k(\eta)}{r^2+s^2})}} 
\bpm 1& \bar f\cr f& |f|^2
\epm,
$$
where $f= e^{2s(h(\xi)-\frac{k(\eta)}{r^2+s^2}) + 2i
r(h(\xi)+\frac{k(\eta)}{r^2+s^2})}$. 
Then 
\begin{align*}
u&= 2is[\ti\pi, a],\\
v&= \frac{1}{|z|^2} (\bar z\ti\pi + z\ti\pi^\perp)b(\eta)(z\ti\pi + \bar z\ti
\pi^\perp).
\end{align*}
The wave map corresponding to $(a, u, v)$ is  
$$\ti s= g_{z,\pi}\ast s= g_{z,\pi}(-1)e^{-ax}(\ti\pi(x,t)-\ti\pi(x,t)^\perp)e^{-ax} g_{z,\pi}(-1).$$
The first column $S$ of $e^{-ax}(\ti\pi(x,t)-\ti\pi(x,t)^\perp)e^{-ax}$ is
$$S(x,t)= \bpm 
-e^{-2imx}\tanh A\cr e^{iB} \sech A\epm,$$
where 
\begin{align*}
A= 2s\left(h\left(\frac{x+t}{2}\right) - \frac{1}{r^2+s^2}\  k\left(\frac{x-t}{2}\right)\right), \\ 
B= 2r\left(h\left(\frac{x+t}{2}\right) +\frac{1}{r^2+s^2}\  k\left(\frac{x-t}{2}\right)\right).
\end{align*}
If $h(x)$ and $k(x)$ tends to zero as $|x|\to \infty$, then   $\lim_{|x|\to \infty} \ti \pi(x,t)= \frac{1}{2}\bpm 1&1\\
1&1\epm$.  In this case the wave maps
$\ti s$ and $S$ have constant boundary condition.   

We claim that if $h'$ and $k'$ are square integrable then the wave map $\ti s$ has finite energy at 
each level $t$.  This can be proved using Theorem \ref{ah}.    The energy at $t=t_0$ is
$$\int_{t=t_0} ||\ti s^{-1}\ti s_x||^2 dx= \int_{t=t_0} \frac{1}{2}||\ti s^{-1}\ti s_\xi + \ti s^{-1} \ti s_\eta||^2 dx.$$
Recall that the trivialization $\ti E$ of the Lax pair of $(a,u,v)$ satisfies the reality condition \eqref{ch}. So $\ti E(x,t,\pm 1)\in U(2)$.  Thus by Theorem \ref{ah}, we have
$$||\ti s^{-1} \ti s_\xi||= 2||h'||, \quad ||\ti s^{-1}\ti s_\eta||= 2||v||.$$

Since $\frac{1}{|z|} (\bar z \ti\pi + z\ti\pi^\perp)\in U(2)$, 
$$|| v(\xi, \eta)||= || k'(\eta)||.$$
The assumption $h', k'$ are in $L^2$ implies that 
$$\int_{t=t_0} ||\ti s^{-1}\ti s_x||^2 dx \ < \ \infty.$$
\eeg

\section{Homoclinic wave maps}\label{ff}

We study the asymptotic behavior of the periodic $k$-soliton wave maps $ s_j$ constructed in Example \ref{bk},
and prove that they are homoclinic. 

 We use the same notations as in Example \ref{bk}.  First we look at the behavior of $ s_j$ as $t\to -\infty$. 
Set 
$$y_j= e^{2\mu_j t}, \quad f_j= e^{2ir_jx}.$$
Since $\mu_j>0$,  $\lim_{t\to-\infty} y_j=0$.  By \eqref{by1}, $q_1(x,t)$ is parallel to
$$\hat q_1= \bpm 1\cr f_1 y_1\epm = \bpm 1\cr 0\epm + f_1y_1\bpm 0\cr 1\epm.$$
So the projection $\pi_1(x,t)$ onto $\C \hat q_1(x,t)$ is
$$\pi_1(x,t)= \bpm 1&0\cr 0&0\epm + y_1\bpm 0& \bar f_1 \cr f_1 &0\epm + O(y_1^2),$$
where $O(y_1^2)$ means terms involving $y_1^n$ with $n\geq 2$.   
By \eqref{bz}, $q_2(x,t)$
 is parallel to 
$$\hat q_2(x,t)= \ti g_1(\bar z_2) \bpm 1\cr f_2y_2\epm.$$
Expand it in $y_1, y_2$ to get  
$$\hat q_2(x,t) = \bpm 1\cr 0\epm + (c_1f_1y_1+ c_2f_2y_2) \bpm 0\cr  1\epm + O(y^2),$$
where $c_1, c_2$ are constants depending on $z_1, z_2$.  
This implies that the projection 
$$\pi_2(x,t )= \bpm 1&0\cr 0&0\epm + y_1\bpm 0&\bar \a_1 \bar f_1 \cr \a_1 f_1 &0\epm 
+ y_2 \bpm 0&\bar \a_2 \bar f_2\cr \a_2 f_2 &0\epm + O(y^2),$$
where $\a_1, \a_2$ are constants depending on $z_1, z_2$.  
Use the formula for $q_j, \pi_j, E_j$ and induction to see that as $t\to -\infty$ we have
\begin{equation}\label{di}
\pi_j(x,t)= \bpm 1&0\cr 0&0\epm + \sum_{r=1}^j  y_r\bpm 0&\bar \a_r \bar f_r \cr \a_r f_r
&0\epm + O(y^2),
\end{equation}
where $\a_r$ are constants depending on $z_1, \ldots, z_j$.  In particular, this proves that
\begin{equation}\label{cb}
\lim_{t\to -\infty} \pi_j(x,t)= \bpm 1&0\cr 0&0\epm.
\end{equation}
Next we use formula \eqref{ca} to compute the asymptotic behavior of $s_j$ as $t\to -\infty$. 
It follows from \eqref{cb} 
that we have
$$\lim_{t\to -\infty} \ti g_r(x,t,-1)^{-1} \ti g_r(x,t,1)=\bpm 1& 0\cr 0 & \frac{1+\bar
z_r}{1+z_r}\epm \ \bpm 1&0\cr 0& \frac{1-z_r}{1-\bar z_r}\epm =
\bpm 1&0\cr 0&-1\epm.$$ By \eqref{ca}, we have
$$\lim_{t\to -\infty}  s_j(x,t)= b_je^{-2ax} \bpm 1&0\cr 0&-1\epm^j c_j^{-1},$$
where $b_j, c_j$ are given by \eqref{ck}.
In particular, we get
\begin{subequations}\label{cm}
\begin{gather}
\lim_{t\to -\infty}  s_{2k}(x,t)= b_{2k}e^{-2ax}c_{2k}^{-1}, \label{cm1}\\
\lim_{t\to-\infty} s_{2k+1}(x,t)= b_{2k+1} e^{-2ax} \bpm 1&0\cr 0&-1\epm
c_{2k+1}^{-1}.\label{cm2}
\end{gather}
\end{subequations}

Both $b_j$ and $c_j$ lies in $U(2)$.  
A direct computation gives
$$g_{z,\pi}(-1) g_{z,\pi}(1)^{-1}= \pi -\pi^\perp.$$
This implies $b_{2k}c_{2k}^{-1}\in SU(2)$.   So 
$$\hat s_{2k}= b_{2k}^{-1} s_{2k} c_{2k}$$
is a wave map from $S^1\times \R^1$ to $SU(2)$.    
Use \eqref{ca} and \eqref{di} and a direct computation to conclude
\begin{equation*} 
\hat s_{2k}(x,t) = e^{-2ax} + e^{-ax} \left(\sum_{j=1}^{2k}  y_j \bpm 0& \bar\b_j\bar f_j\cr - \b_jf_j
&0\epm \right) e^{-ax} + O(y^2) 
\end{equation*}
as $t\to-\infty$, 
where $\b_1, \ldots, \b_{2k}$ are real constants. Recall that 
$$s_0(x,t)= e^{-2ax} = \diag(e^{-2imx}, \ e^{2imx}).$$
So we have 
$$s_0^{-1}(x,t) y_j e^{-ax} \bpm 0& \bar \b_j\bar f_j\cr -\b_j f_j
&0\epm e^{-ax} = e^{2\mu_j t}\bpm 0 & \bar\b_j e^{2i(m-r_j)x}\cr -\b_j e^{-2i(m-r_j) x} & 0\epm,$$
which is an unstable mode $p^+_{-2m, 2r_j}$  at the stationary solution
$s_0$ with eigenvalue
$2\sqrt{m^2-r_j^2}=2\mu_j$ as given in Corollary \ref{cza}.  In other words, we have shown 
\begin{equation} \label{cc}
\lim_{t\to-\infty} \hat s_{2k}(x,t)- s_0(x,t) - s_0(x,t)\sum_{j=1}^{2k} 
p_{-2m,2r_j}^+(x,t) =0.
\end{equation}

To compute the asymptotic behavior of $s_j(x,t)$ as $t\to  \infty$, we set
$$h_j= e^{-2ir_j x},\quad \rho_j= e^{-2\mu_j t}.$$
Since $\mu_j>0$, $\lim_{t\to \infty} \rho_j=0$.  
A similar computation  implies that  
\ben
\item $q_1(x,t)$ is parallel to  $\bpm h_1\rho_1 \cr 1\epm$,  

\item 
\begin{equation}\label{cd}
\pi_j(x,t)=  \bpm 0&0\cr 0&1\epm  + \sum_{n=1}^j  \rho_n \bpm 0& \bar \b'_n\bar h_n\cr \b'_n h_n
&0\epm +O(\rho^2)
 \end{equation}
for some constants $\b'_1, \ldots, \b'_j$ depending on $z_1, \cdots, z_j$,

\item $\lim_{t\to \infty} \ti g_j(x,t,-1)^{-1} \ti g_j(x,t, 1)=\bpm -1&0\cr 0&1\epm$,
\item 
\begin{subequations}\label{cn}
\begin{gather}
\lim_{t\to \infty} s_{2k}(x,t)= b_{2k} e^{-2ax} c_{2k}^{-1}, \label{cn1}\\
\lim_{t\to \infty}s_{2k+1}(x,t)= b_{2k+1} e^{-2ax}\bpm -1 &0\cr 0&1\epm c_{2k+1}^{-1}.\label{cn2}
\end{gather}
\end{subequations}

\item
\begin{equation*}
\hat s_{2k}(x,t)=e^{-2ax} +e^{-ax} \sum_{j=1}^{2r}  \rho_j 
\bpm 0& \bar \e_j\bar h_j\cr -\e_j h_j &0 \epm e^{-ax} + O(\rho^2), \quad {\rm as\ } t\to \infty
\end{equation*}
for some constants $\e_1, \ldots, \e_j$ depending on $z_1, \cdots, z_j$. 
\item As $t\to \infty$ we  have
\begin{equation}\label{ce}
\lim_{t\to \infty} \hat s_{2k}(x,t)- s_0(x,t) - s_0(x,t)\sum_{j=1}^{2k} p^-_{-2m, 2r_j}(x,t) =0
\end{equation}
for some linear stable mode $p^-_{-2m, 2r_j}$.  
\een

As a consequence of \eqref{cm2} and \eqref{cn2}, we get
$$\lim_{t\to -\infty} s_{2k+1}(x,t)= -\lim_{t\to \infty} s_{2k+1}(x,t)= b_{2k+1}e^{-2ax} \bpm 1&0\cr
0& -1\epm c_{2k+1}^{-1}.$$
So $s_{2k+1}$ is a hecteroclinic wave map from $S^1\times \R^1$ to $SU(2)$.  

Formulas \eqref{cc}
and \eqref{ce} imply that 
$s_{2k}, \hat s_{2k}$ are homoclinic wave map from $S^1\times \R^1$ to $SU(2)$.  So we have

\begin{thm} \label{bb}
Let $2m$ be an integer, $a=\diag(im, -im)$, and $s_0(x,t)= e^{-2ax}$.  Let $\pi$ denote the
projection of $\C^2$ onto $\C \bpm 1\cr 1\epm$, and 
$$z_j= e^{i\o_j}= \frac{r_j}{m} + i\ \frac{\sqrt{m^2-r_j^2}}{m} = \frac{r_j+i \mu_j}{m},$$
where $2r_1, \ldots, 2r_{2k}$ are integers and $|r_j|<m$.  Let 
$$s_{2k}= g_{z_{2k}, \pi}\bu(\cdots \bu (g_{z_1, \pi}\bu s_0)\cdots)$$
 denote the wave
map obtained by applying $2k$ B\"acklund transformations.
Then $s_{2k}$ is a homoclinic wave map from $S^1\times \R^1$ to $SU(2)$.  Moreover,  there exist
constants $b_{2k}, c_{2k}\in U(2)$ such that
$\hat s_{2k}= b_{2k}^{-1}s_{2k} c_{2k}$ satisfies  $\lim_{|t|\to \infty} s_{2k}(x,t)= s_0(x,t)$ and
\eqref{cc} and \eqref{ce}.
\end{thm}

\section{wave maps into $S^2$}\label{fg}

We describe a constraint condition for solutions the $-1$ equation associated to $SU(2)=S^3$ so that the
corresponding wave maps into $S^3$ actually lies in $S^2$.    

Recall that we identify $S^3\subset \C^2=\R^4$ as $SU(2)$ via 
$$\bpm z\cr w\epm \mapsto \bpm z& -\bar w\cr w& \bar z\epm.$$ 
The intersection of $S^3$ with the linear hyperplane defined by {\rm Re}$(w)=0$ is
\begin{equation}\label{cv}
M=\left\{\bpm z& is\cr is & \bar z\epm\ \bigg|\ \   z\in \C, \ s\in \R,\  |z|^2+ s^2=1\right\},
\end{equation}
which is a totally geodesic $2$-sphere in
$SU(2)=S^3$. If $y\in SU(2)$, then $y\in M$ if and only if $y^t=y$.   

It is well-known that if
$s$ is a wave map into a Riemannian manifold $N$ and the image of $s$ lies in a totally geodesic
submanifold $M$ of $N$, then $s$ is a wave map into $M$.  Hence if a wave map into $S^3$
with its image lies in a totally geodesic $S^2$, then this is also a wave map into
$S^2$. 

In this section, we give a constraint condition on solutions of the $-1$ flow equation so that the corresponding
wave maps lie in the totally geodesic $S^2$ in $SU(2)=S^3$.  As a consequence, we also see that the SGE is a
subequation of the constrained $-1$ flow equation.       

In order to explain how to get wave maps into $S^2$, we need to view $S^2$ as the
symmetric space $SU(2)/SO(2)$.  Let
$\tau, \sigma:SL(2,\C)\to SL(2,\C)$ denote the  maps defined by 
$$\tau(g)= (g^*)^{-1}, \quad \sigma(g)=(g^t)^{-1}.$$
Then $\sigma$ and $\tau$ are group homomorphisms, $\sigma^2=\tau^2={\rm id}$, $\sigma
\tau=\tau\sigma$, and the differentials
 at the identity matrix are
$$\sigma_\ast(\xi)= -\xi^t, \quad \tau_\ast(\xi)=-\xi^*.$$
Note that $\sigma(SU(2))\subset SU(2)$.  Hence $\sigma |SU(2)$ is an involution of
$SU(2)$. The fixed point set of $\sigma$ in $SU(2)$ is $SO(2)$, and $S^2$ is diffeomorphic
to $SU(2)/SO(2)$. A direct computation shows that
$$M_\sigma=\{g\sigma(g)^{-1}\n g\in SU(2)\}$$
is the totally geodesic $M=S^2$ given in \eqref{cv}.  

Since $\sigma_*^2={\rm id}$, we have:
$$su(2)= \ck+\cp,$$  where $\ck$ and $\cp$ are  eigenspaces of
$\sigma_\ast$ on
$su(2)$ with eigenvalues $1, -1$ respectively.  In fact,   
\begin{equation}\label{er}
\ck= so(2), \quad \cp= \left\{ i\bpm x& y\cr y& x\epm\ \bigg|\  x, y\in\R\right\}.
\end{equation}
Since $\sigma_*$ is a Lie algebra automorphism,  $\sigma_\ast([\eta_1, \eta_2])=
[\sigma_*(\eta_1),
\sigma_*(\eta_2)]$.  So we have
$$[\ck, \ck]\subset \ck, \quad [\ck, \cp]\subset \cp, \quad [\cp, \cp]\subset \ck.$$

\bprop  \label{cp}
Let $(a,u,v)$ be a solution of the $-1$ flow equation associated to $SU(2)$, $\o_\l= (a\l + u) \ d\xi +
\l^{-1} \ d\eta$ its Lax pair, and $E$ the trivialization of $\o_\l$.  Then the following statements are
equivalent:
\ben 
\item $a, v\in \cp$ and $u\in \ck$.
\item  $\o_\l$ satisfies 
\begin{equation} \label{cr}
\o_\l= - \o_{\bar\l}^*, \quad \o_\l= -\o_{-\l}^t.
\end{equation}
\item The trivialization $E(x,t,\l)$ of $\o_\l$ satisfies 
\begin{equation}\label{cs}
E(x,t,\bar\l)^*E(x,t,\l)= {\rm I\/}, \quad E(x,t,\l)E(x,t,-\l)^t={\rm I\/},
\end{equation}
\een
\eprop

\begin{proof}
It is easy to see that (1) and (2) are equivalent and (3) implies (2).  To prove (2) implies (3), we let
$F(x,t,\l)= (E(x,t,-\l)^t)^{-1}$.  Then 
$$F^{-1} dF= -\o_{-\l}^t= \o_\l, \quad F(0,0,\l)={\rm I\/}.$$ Hence $E=F$, i.e., 
$E(x,t,\l)E(x,t,-\l)^t= {\rm I\/}$. 
We have shown before that $\o_\l=-\o_{\bar \l}^*$ implies that
$E(x,t,\bar\l)^*E(x,t,\l)={\rm I\/}$.  
\end{proof}

\bprop \label{dy}
Let $(a,u,v)$ be a solution of the $-1$ flow equation associated to $SU(2)$, and $s$ the wave map
constructed from $(a,u,v)$ in Theorem \ref{ah}.  If $a, v\in \cp$ and $u\in \ck$, then $s$ is a wave
map into $S^2$. 
\end{prop}

\begin{proof} 
By Proposition \ref{cp}, the trivialization $E$ of the Lax pair corresponding to $(a,u,v)$ satisfies the reality
condition \eqref{cs}.
This implies that
$E(x,t,r)\in SU(2)$ for all
$r\in \R$ and $E(x,t,1)^{-1}= E(x,t,-1)^t$.  The wave map constructed in Theorem \ref{ah} is 
$$s(x,t)= E(x,t,-1)E(x,t,1)^{-1}= E(x,t,-1) E(x,t,-1)^t.$$ 
But $y\in SU(2)$ lies in 
$S^2$ given in \eqref{cv} if and only if $y^t=y$.   So $s(x,t)$ lies in
$S^2$.  Because $S^2$ is totally geodesic in $SU(2)$, $s$ is a wave map into $S^2$. 
\end{proof}

The following Proposition was proved in \cite{TerUhl00a}.  
\bprop \label{cq}
Let $(a, u,v)$ be a solution of the $-1$ flow equation associated to $SU(2)$ and $a, v\in
\cp$ and $u\in\ck$, and $s$ the wave map corresponding to $(a,u,v)$ in Theorem
\ref{ah}.  Then:
\ben
\item[(1)] 
If $z=i\mu$ is pure imaginary and $\bar \pi=\pi$, then $(a, \ti u, \ti v)= g_{i\mu,\pi}\bu
(a,u,v)$ is again a solution of the $-1$ flow with $\ti u\in \ck$, $\ti v\in \cp$,   and
$g_{i\mu,\pi}\bu s$ is a wave map into $S^2$.   
\item[(2)] 
If $\bar \pi=\pi$, then
$(a, \hat u, \hat v)= g_{z,\pi}\bu (g_{-\bar z, \pi}\bu (a,u,v))$ is a solution of the $-1$ flow
with $\hat u\in \ck$, $\hat v\in \cp$, and $g_{z,\pi}\bu (g_{-\bar z, \pi}\bu s)$ is a wave
map into $S^2$.  
\een
\eprop

\beg  {\bf SGE and wave maps into $S^2$} (\cite{Poh76, ShaStr96, TerUhl00a}).  
\par

Let  $su(2)=\ck+\cp$, where $\ck$ and $\cp$ are given by \eqref{er}.  Let
$$a=\diag(i, -i), \quad u=\bpm 0 & \frac{q_x}{2}\cr -\frac{q_x}{2}& 0\epm, \quad
 v=-\frac{i}{4}\bpm \cos q& \sin q\cr \sin q& -\cos q\epm.$$
Note $a, v\in \cp$ and $u\in \ck$.
A direct computation implies that:
\ben
\item[(1)] $(a, u, v)$ is a solution of the $-1$ flow equation associated
to $SU(2)$ if and only if $q$ is a solution of the sine-Gordon equation (SGE):
\begin{equation}\label{cy}
q_{xt}= \sin q.
\end{equation}
Hence solutions of the SGE give rise to wave maps into $S^2$. 
\item[(2)] Let  $(a, \ti u, \ti v)= g_{is,\pi}\bu (a, u,v)$.   If $\bar\pi=\pi$, then 
$$\ti u=\bpm 0 & \frac{\ti q_x}{2}\cr -\frac{\ti q_x}{2}& 0\epm, \quad
\ti v=-\frac{i}{4}\bpm \cos \ti q& \sin \ti q\cr \sin \ti q& -\cos \ti q\epm$$
for some $\ti q$.  
So $\ti q$ is again a solution of the SGE.  Let $g_{is,\pi}\bu q$  denote $\ti q$. 
\item[(3)] $g_{is, \pi}\bu 0$ is the traveling wave
solution of the SGE:
$$q(x,t)= 4\tan^{-1}(e^{sx+\frac{t}{ s}}),$$
and $g_{e^{i\o}, \pi}\bu (g_{-e^{-i\o}, \pi}\bu 0)$ is the breather solution of the SGE:
$$4\tan^{-1}\left(\frac{\sin\o\sin (t\cos\o )}{ \cos\o
\cosh(x\sin\o)}\right).$$
\een
This solution is periodic in time, and hence called a breather.  But wave map equation is
invariant if we exchange the space and time variables.  So a breather solution can also be viewed
as periodic in space instead.  
\eeg

Since wave maps into $S^2$ is a special cases of wave maps into $SU(2)$, as a consequence
of Theorem \ref{bb} and Proposition \ref{cq} we get the following Theorem:

\begin{thm}\label{de}
We use the same assumption as in Theorem \ref{bb}.  If
$z_{2j}= -\bar z_{2j-1}$, i.e., $r_{2j}= -r_{2j-1}$ and $\mu_{2j}= \mu_{2j-1}$ for all $1\leq
j\leq k$, then $s_{2k}$ is a homoclinic $2k$-soliton wave map from $S^1\times\R^1$ into $S^2$. 
\end{thm}

When $k=1$,  the above Theorem was proved by Shatah and Strauss in \cite{ShaStr96}.

\section{Wave maps into compact symmetric spaces\/}\label{fh}

In section \ref{fg}, we embed $S^2$ as a totally geodesic submanifold of $SU(2)=S^3$. By viewing
$S^2$ as the symmetric space $SU(2)/SO(2)$, we give conditions on solutions of the $-1$
flow equation whose image lies in $S^2$. In fact, this same method works for any compact
symmetric space.  In particular, we apply this method to construct homoclinic periodic $2k$-soliton wave maps into
$\C P^n$ and $4k$-soliton wave maps into $S^n$.  

First we give a short review of symmetric spaces.  Let $G$ be a complex
semi-simple Lie group, and $\tau$ and $\sigma$ involutions of $G$ such that
\ben
\item[(i)] 
the differential
$\tau_\ast= d\tau_e$ and $\sigma_\ast= d\sigma_e$ at the identity $e$ are conjugate linear
and complex linear Lie algebra involution on $\cg$ respectively, i.e.,
$\tau_*(\a \xi)=\bar \a \tau_*(\xi)$ and $ \sigma_*(\a \xi)= \a \sigma_*(\xi)$
for all $\a\in \C$ and $\xi\in \cg$,
\item[(ii)]
$\sigma\tau= \tau\sigma$.
\een
Let $U$ denote the fixed point set of $\tau$ in $G$.  Such $U$ is called a {\it real  form of $G$\/}.  Since
$\tau$ and
$\sigma$ commute,
$\sigma(U)\subset U$.  Let $K$ denote the fixed point set of $\sigma$ in $U$, and $\cp$ the
$-1$ eigenspace of $\s_\ast$ on $\cu$.  Then
$U/K$ is a symmetric space, and 
$\cu=\ck + \cp$ satisfying 
$$[\ck, \ck]\subset \ck, \quad [\ck, \cp]\subset \cp, \quad [\cp, \cp]\subset \ck.$$

Let $D$ be a domain in $\C$ that is invariant under complex conjuation.  We say that $g:D\to G$ and
$\xi:D\to \cg$ satisfy
 the  {\it $U$-reality condition\/} if
\begin{equation}\label{dc}
\tau(g(\bar\l))=g(\l), \quad \tau_*(\xi(\bar\l))=\xi(\l),
\end{equation} 
and satisfy the {\it $U/K$-reality condition\/} if
\begin{equation}\label{dd}
\bca
\tau(g(\bar\l))=g(\l), \quad \sigma(g(-\l))= g(\l),& \\
\tau_*(\xi(\bar\l))= \xi(\l), \quad \sigma_*(\xi(-\l))= \xi(\l)
\eca
\end{equation}
respectively.  

Let $\ast$ denote the $U$-action on $U$ defined by $g\ast h= gh\sigma(g)^{-1}$.  Then the stablizer at 
the identity $e$ is $K$.  So the orbit at $e$, 
\begin{equation}\label{df}
M_\sigma=\{g\sigma(g)^{-1}\n g\in U\},
\end{equation}
 is diffeomorphic to $U/K$. It is known that $M_\sigma$ is isometric to the symmetric space $U/K$:

\begin{prop}
The $U$-orbit $M_\sigma= U\ast e$ is totally geodesic submanifold of $U$ and is an isometric
embedding of the symmetric space $U/K$ into $U$.  
\end{prop}

The embedding of $U/K$ given in the above Proposition is called the {\it Cartan
embedding\/}.  The embedding of $S^2$ in $SU(2)=S^3$ given by \eqref{cv} is the Cartan 
embedding of $M_\sigma$ in $SU(2)$, where $\frac{SU(2)}{SO(2)}=S^2$ is the symmetric space given by
$\tau(g)= (g^*)^{-1}$ and  $\sigma(g)= (g^t)^{-1}$. 

{\it The $-1$ flow equation associated to $U$\/} is the equation \eqref{ac} for
$(a,u,v):\R^2\to
\prod_{i=1}^3 \cu$, and has a Lax pair $\o_\l= (a\l+ u)\ dx+ \l^{-1} v \ dt$.  Our computations and results for
$SU(2)$ and $S^2=\frac{SU(2)}{SO(2)}$ in previous sections work for any compact Lie group $U$ and
symmetric space $U/K$.  For example, the following can be proved in a similar manner:
\ben

\item The Lax pair $\o_\l$ of the $-1$ flow equation associated to $U$ satisfies the $U$-reality condition
\eqref{dc}, i.e., 
$$\tau_*(\o_{\bar\l})= \o_\l.$$

\item Theorem \ref{ah} holds by replacing $SU(n)$ by $U$. In other words, we have a
correspondence between solutions of the $-1$ flow equation associated to $U$ and wave maps into $U$.

\item Since any compact Lie group $U$ can be embedded as a subgroup of $SU(N)$ for some $N$, to
construct explicit solutions of the $-1$ flow equation, we only need to find product of simple elements
(i.e.,  of the form \eqref{do}) that satisfies the $U$-reality condition.  

\item The linearization of the wave map equation from $S^1\times \R^1$ into $U$ at the stationary
wave map
$s_0(x,t)= e^{ax}$ with $e^{2\pi a}={\rm I\/}$ is 
$$\xi_{tt}=\xi_{xx} + [a, \xi_x].$$
Its stable and unstable modes can be computed using roots of $\cu$.  

\item  Propositions \ref{cp} and \ref{dy} hold if we replace $SU(2)$ and 
$\frac{SU(2)}{SO(2)}$ by $U$
and
$\frac{U}{K}$.  The proofs are similar.  

\item There are analogous Theorem \ref{bb} and Theorem \ref{de} for $U$ and $U/K$.  
\een

Next we give two examples:

\beg\label{dn} {\bf Wave maps from $S^1\times \R^1$ to $\C P^{n-1}$}
\par

Let $G= SL(n,\C)$, $J=\diag(1, \ldots, 1, -1)$, and $\tau, \sigma: G\to G$ defined by
$$\tau(y)= (y^*)^{-1}, \quad \sigma(y)=JyJ^{-1}.$$ 
A direct computation shows that both $\tau$ and $\sigma$ are group homomorphisms, 
$\tau^2=\sigma^2={\rm Id}$, and
$\tau\sigma=\sigma\tau$.  The fixed point set of $\tau$ in $G$ is $U=SU(n)$, and the fixed point set of
$\sigma$ in $U$ is
$S(U(1)\times U(n-1))$.  The $\pm 1$ eigenspace of $\sigma_\ast$ in $\cu$ is 
\begin{align*}
\ck&= \left\{\bpm \xi &0\cr 0& c\epm \bigg|\  \xi\in u(n-1), \ c\in \C \ {\rm pure\
imaginary\/}, \ {\rm tr\/}(\xi)+ c=0\right\}\cr
\cp&=\left\{\bpm 0& \bar v^t\cr -v &0\epm\bigg| \ v \in  \cm_{1\times (n-1)}(\C)\right\}.
\end{align*}
Here $\cm_{k\times j}(K)$ is the space of $k\times j$ matrices with entries in $K$.  
 The symmetric space corresponding to $\tau, \sigma$ is $\frac{SU(n)}{S(U(n-1)\times
O(1))} = \C P^{n-1}$. 

Let $D\in gl(n-1,\C)$, $v\in \cm_{(n-1)\times 1}(\C)$, $v\in \cm_{1\times (n-1)}(\C)$,
$c\in \C$.  Then
$g=\bpm D& b\cr v & c\epm \in SU(n)$ if and only if 
$$\bca D\bar D^t+b\bar b^t= {\rm I\/}, &\\
v\bar D^t + c\bar b^t=0,&\\
||v||^2 + |c|^2=1.
\eca $$
So 
$$g\s(g)^{-1}= \bpm I- 2b\bar b^t& 2b\bar c\cr - 2c \bar b^t& 2|c|^2-1\epm.$$
The map from the symmetric space $M_\sigma$ to $\C P^{n-1}$ (the space of complex linear lines in $\C^n$)
defined by
$$ \bpm I- 2b\bar b^t& 2b\bar c\cr - 2c \bar b^t& 2|c|^2-1\epm\mapsto  \C\bpm 2 b\bar c\cr 2
|c|^2-1\epm$$  is an isometry from $M_\sigma$ to $\C P^{n-1}$.  

 In order to construct explicit wave maps into $\C P^{n-1}$, we look for a product of simple elements, 
$$g=g_{z_1, \pi_1}\cdots g_{z_r, \pi_r},$$
 that satisfies the extra reality condition 
$$g(-\l)= \sigma(g(\l)).$$
It can be checked that the simple element $g_{z,\pi}$ does not satisfy this extra reality
condition.  But we can find product of two simple elements do.   Let
$\pi$ denote the Hermitian projection onto
$\C\bpm w\cr c\epm$, where
$w\in
\C^{n-1}$ and $c\in \C$ so that $||w||^2= |c|^2=1$.  A direct computation implies
\begin{subequations}\label{dp}
\begin{gather}
\pi= \frac{1}{2}\bpm w\bar w^t & w\bar c\cr c\bar w^t & |c|^2\epm, \label{dp1}\\
\sigma(\pi)=  \frac{1}{2}\bpm  w\bar w^t & -w\bar c\cr -c\bar w^t & |c|^2\epm = {\rm projection\
onto\ } \C\bpm w\cr -c\epm.\label{dp2}
\end{gather}
\end{subequations}
Since $\bpm w\cr c\epm$ and $\bpm w\cr -c\epm$ are perpendicular with respect to the Hermitian
inner product, we have
\begin{equation}\label{dq}
\pi\sigma(\pi)= \sigma(\pi)\pi=0.
\end{equation}
So
$g_{z,\pi}$ and $g_{-z, \sigma(\pi)}$ commute.  Let 
\begin{equation}\label{dr}
h_{z,\pi}= g_{z,\pi}g_{-z, \sigma(\pi)}.
\end{equation}
Note
$$\sigma(\pi^\perp) = \sigma({\rm I\/} - \pi)= {\rm I\/} - \sigma(\pi)=(\sigma(\pi))^\perp.$$
Use \eqref{dq}, $\sigma(\pi^\perp)= \sigma(\pi)^\perp$ and a  direction computation to prove that
$h_{z,\pi}$ satisfies the
$\frac{SU(n)}{S(U(n-1)\times U(1))}$-reality condition:
\begin{equation}\label{eb}
h(\bar\l)^*h(\l)= {\rm I\/}, \quad h(-\l) = \sigma(h(\l)).
\end{equation}
We apply B\"acklund transformations given by these elements to the stationary wave maps (closed
geodesics) to construct homoclinic wave maps from $S^1\times \R$ to $\C P^{n-1}$. 

Let $a_0\in \cp$ so that $e^{2\pi a_0}={\rm I}$, $m$ an integer, and $a= ma_0$.  Since $a_0\in
su(n)$, there exists $A\in SU(n)$ and $C=\diag(ic_1, \ldots, ic_n)$ so that 
$a_0= ACA^{-1}$.  Because $e^{2\pi a_0}={\rm I\/}$, all $c_j$ must be integers.   
Note $(a,0,a)$ is a solution of the $-1$ flow equation associated to $SU(n)$ with $a\in\cp$, and the
corresponding wave map constructed in Theorem \ref{ah} is 
$$f_0(x,t)= e^{-2ax}= Ae^{-2Cx}A^{-1}\in M_\sigma = \C P^{n-1}.$$
Choose $z_1, \ldots, z_j\in \C$ so that
$$z_j=e^{i\o_j}=\frac{r_j+i\mu_j}{m}$$
as in Example \ref{bk}.  Let 
\begin{align*}
 h_{z,\pi}\bu (a,u,v)&= g_{z,\pi}\bu (g_{-z, \sigma(\pi)}\bu (a,u,v)), \\
(a, u_j, v_j)  &= h_{z_j, \pi}\bu(\cdots \bu (h_{z_1, \pi}\bu (a,0,a))\cdots ), \\
f_j &= h_{z_j, \pi}\bu(\cdots \bu (h_{z_1, \pi}\bu f_0)\cdots ).
\end{align*}
The computation in Example \ref{bk} and the proof of Theorem \ref{bb} implies that 
$f_k$ is a homoclinic $2k$-soliton wave map from $S^1\times \R^1$ to $SU(n)$.  But each $h_{z_j,\pi_j}$
satisfies the reality conditions \eqref{eb}. Use a proof similar to that of Proposition \ref{cq}(2)
to see  that $v_j\in \cp$ and $u_j\in \ck$ and the wave map corresponding to
$(a, u_j,v_j)$ lies in $M_\sigma= \C P^{n-1}$.  So $f_j$ is a wave map from $S^1\times \R$
to $SU(n)$ whose image lies in $\C P^{n-1}$.  By Theorem \ref{bb}, $f_{k}$ is a homoclinic wave map
from $S^1\times \R$ to $\C P^{n-1}$.  

\eeg

\beg \label{dg} {\bf Wave maps from $S^1\times \R^1$ to $\R P^{n-1}$ and $S^{n-1}$}
\par

Let $G=SO(n,\C)$, 
$$\tau(g)= \bar g, \quad \sigma(g)= JgJ^{-1},$$ where $J= \diag(1, \ldots, 1, -1)$. 
It can be checked easily that both $\tau$ and $\sigma$ are group homomorphisms of $G$, 
$\tau^2=\sigma^2={\rm Id}$ and
$\tau\sigma=\sigma\tau$.  The fixed point set of $\tau$ in $G$ is $U=SO(n)$, and the fixed point set of
$\sigma$ in $U$ is
$S(O(1)\times O(n-1))$.  The $\pm 1$ eigenspace of $\sigma_\ast$ in $\cu$ is 
\begin{align*}
\ck&= \left\{\bpm 0 &0\cr 0& \xi\epm \bigg| \ \xi\in so(n-1)\right\}\cr
\cp&=\left\{\bpm 0& v^t\cr -v &0\epm\bigg|\ v \in \cm_{1\times (n-1)}(\R)\right\}.
\end{align*}
 The symmetric space corresponding to $\tau, \sigma$ is $\frac{SO(n)}{S(O(n)\times
O(1))}$, which is $\R P^{n-1}$.  Let
$g=\bpm D&b\cr v&c\epm\in SO(n)$ with $D\in gl(n-1,\R)$, $b\in
\cm_{(n-1)\times 1}*(\R)$, $v\in \cm_{1\times(n-1)}(\R)$, and $c\in \R$. A direct computation gives
$$g\sigma(g)^{-1}= \bpm {\rm I\/}- 2bb^t & 2bc\cr -2cb^t & 2c^2-1
\epm.$$
So the Cartan embedding is 
\begin{align*}
M_\sigma & = \{g\sigma(g)^{-1}= gJgJ^{-1}\n g\in SO(n)\} \cr
&=\left\{\bpm {\rm I\/} - 2bb^t& 2bc\cr -2cb^t & 2c^2-1\epm\bigg|\,\,
b\in \cm_{(n-1)\times 1}(\R), \ c\in \R,\   c^2+||b||^2=1\right\}.
\end{align*}
The map from 
$$S^{n-1}=\left\{\bpm b\cr c\epm \bigg|\  c\in \R, \ b\in \R^{n-1},\  |c|^2+
||b||^2=1\right\}$$ to
$M_\sigma$ defined by $$\bpm b\cr c\epm \ \mapsto \bpm {\rm I\/} - 2bb^t& 2bc\cr -2cb^t
& 2c^2-1\epm$$ is a double covering.  This shows that $M_\sigma$ is isometric to $\R
P^{n-1}$.

To construct explicit wave maps into $S^{n-1}$ from a stationary wave map into $S^{n-1}$,
we need to find rational maps
$g$ from $S^2$ to $SO(n,\C)$ that satisfying $g(\infty)={\rm I\/}$ and the
$\frac{SO(n)}{S(O(n-1)\times O(1))}$-reality conditions:
\begin{equation}\label{ds}
\overline{g(\bar\l)} g(\l)= {\rm I\/}, \quad g(-\l)= Jg(\l) J^{-1}.
\end{equation}
This is equivalent to find rational maps from $S^2$ to $GL(n,\C)$ that satisfies
\begin{equation}\label{dt}
g(\bar \l)^*g(\l)={\rm I\/}, \quad \overline{g(\bar\l)} g(\l)= {\rm I\/}, \quad g(-\l)= Jg(\l) J^{-1}.
\end{equation}
We have seen in Example \ref{dn} that $h_{z,\pi}$ defined by \eqref{dr} satisfies the first and third
conditions of  \eqref{dt}, but in general it does not satisfies the second condition.  However, let $w\in
\R^{n-1}$ and $b\in \R$ so that $||w||= |b|=1$, and $\pi$ the projection onto $\C\bpm w\cr ib\epm$. 
A direct computation shows that
$$\pi= \bpm ww^t& -ibw\cr ibw^t&1\epm, \quad \sigma(\pi)= \bar \pi.$$
By \eqref{dq}, we have
$$\pi \bar\pi=\bar\pi\pi=0.$$
It is easy to check that
$$\phi_{z,\pi}= h_{z,\pi} h_{-\bar z, \pi}= g_{z,\pi}g_{-z, \bar\pi} g_{-\bar z, \pi} g_{\bar z, \bar \pi}$$
satisfies all conditions in \eqref{dt}.  Note $\phi$ has four simple poles, $z, -z, \bar z, -\bar z$.  

Let $a_0\in \cp$ so that $e^{2\pi a_0}= {\rm I\/}$.  Then $(a, 0,a)$ is a solution of the $-1$ flow
equation associated to
$SO(n)$ whose corresponding wave map is the stationary wave map into $M_\sigma$: 
$s_0(x,t)=e^{-2a_0x}$.  Let
$m$ an integer, and 
$$a=ma_0.$$  
    Let 
$$z_j= e^{i\o_j}= \frac{r_j+ i\mu_j}{m},$$ where $r_j$ is an integer. Let 
$$\phi\bu(a,u,v) = h_{z,\pi}\bu (h_{-\bar z, \pi}\bu (a,u,v)).$$
 Use similar reasoning as in Example \ref{dn} to see that wave map 
$$s_j= \phi_{z_{j},\pi} \bu (\cdots (\phi_{z_1, \pi}\bu s_0)\cdots )$$
is a homoclinic $4j$-soliton wave map from $S^1\times \R^1$ to $\R P^{n-1}$. 
Moreover, the last column of $s_j$ gives a homoclinic wave map from $S^1\times \R^1$ to
$S^{n-1}$.   
\eeg

\section{Wave maps from $\R^{1,1}$ to $SL(2,\R)$}\label{fi}

 It is known that the Cauchy problem for wave map equation from $\R^{1,1}$ to any complete Riemannian
manifold has long time existence (\cite{Gu80}). We will show that this is no longer true when the target
manifold is the pseudo-Riemannian manifold $SL(2,\R)$.  In fact, we use B\"acklund
transformations to construct smooth initial data with finite energy and constant boundary condition at $\pm \infty$  so that the  Cauchy problem for wave maps into $SL(2,\R)$ 
\ben
\item has long time existence, or
\item develops singularities in finite time. 
\een

First note that Theorem \ref{ah} holds if we
replace $SU(n)$ by any group  $G$.  But B\"acklund transformations for the $-1$ flow equation associated to $SL(2,\R)$ is
different from the $SU(n)$ case.  
 Let $\a_1, \a_2\in \C$, and $\pi$ a linear projection of $\C^n$ (i.e., $\pi$ is complex linear and
$\pi^2=\pi$).  Let
\begin{equation}\label{eg}
h_{\a_1, \a_2, \pi}(\l)= {\rm I\/} + \frac{\a_1-\a_2}{\l-\a_1}\ \pi',
\end{equation}
where $\pi'={\rm I\/}-\pi$.   Then 
$$h_{\a_1, \a_2, \pi}^{-1}(\l) = {\rm I\/} + \frac{\a_2-\a_1}{\l-\a_2}\ \pi'.$$
B\"acklund transformations are given as follows  (cf. \cite{TerUhl00a}):

\bthm  \label{ef}
Let $(a,u,v)$ be a smooth solution of the $-1$ flow equation associated to $SL(n,\C)$, $\o_\l= (a\l+ u)\
d\xi + v\l^{-1} d\eta$ its Lax pair, and $E$ the trivialization of $\o_\l$.  Let $\pi$ be the projection of
$\C^n$, $V_1={\rm Im\/}(\pi)$, and  $V_2={\rm Ker\/}(\pi)$.  Set
$$\ti V_i(\xi,\eta)= E(\xi,\eta, \a_i)^{-1}(V_i), \quad {\rm for\ } i=1, 2.$$
Suppose $\ti V_1(\xi,\eta)\cap \ti V_2(\xi,\eta)=\{0\}$ for $(\xi,\eta)$ lies in an open subset $\co$ of
$\R^2$.  Let $\ti \pi(\xi,\eta)$ denote the linear projection onto $\ti V_1(\xi,\eta)$ along $\ti V_2(\xi,
\eta)$, and
\begin{align*}
\ti u &= u+ (\a_1-\a_2)[a, \ti \pi], \\
\ti v&= \left({\rm I\/}-\frac{\a_1-\a_2}{\a_1}\ \ti \pi\right) v  \left({\rm I\/}-\frac{\a_1-\a_2}{\a_1}\
\ti \pi\right)^{-1} = (-\a_1\ti\pi + \a_2\ti\pi')v (-\a_1^{-1}\ti\pi + \a_2^{-1} \ti \pi').
\end{align*}
Then $(a,\ti u, \ti v)$ is a smooth solution of the $-1$ flow equation associated to $SL(n,\C)$
defined on $\co$.
Moreover, the trivialization of the Lax pair of $(a, \ti u, \ti v)$ is 
$$\ti E(\xi, \eta, \l) =h_{\a_1, \a_2, \pi}(\l) E(\xi, \eta, \l) h_{\a_1, \a_2, \ti \pi(x,t)}(\l)^{-1}.$$
\ethm

We use $h_{\a_1, \a_2,\pi}\bu (a, u, v)$ to denote $(a, \ti u, \ti v)$.  
   
If $(a,u,v)$ is a solution of the $-1$ flow equation associated to $SL(n,\R)$, then $h_{\a_1, \a_2,
\pi}\bu (a, u,v)$ is also a solution of the $-1$ flow equation associated to $SL(n,\R)$ provided that
$\a_1, \a_2\in \R$ and $\bar\pi=\pi$.  

\beg 
Let $a=\diag(1, -1)$. Then $(a, 0, a)$ is a solution of the $-1$ flow equation associated to $SL(2,\R)$
and the corresponding wave map is $s(x,t)= e^{-2ax}$.  Let
$\a_1, \a_2\in \R$, $v_1=\bpm c_1\cr c_2\epm\in \R^2$, $v_2= \bpm d_1\cr d_2\epm\in \R^2$,
and
$\pi$ the projection of $\C^2$ onto $\C v_1$ along $\C v_2$.  Let 
$$(a, \ti u, \ti v)= h_{\a_1, \a_2,\pi}\bu (a, 0,a).$$
Use Theorem \ref{ef} and a direct computation to get 
$$\ti \pi= \frac{1}{c_1 d_2 e^A- c_2d_1 e^{-A}}\ \bpm c_1 d_2 e^A& -c_1d_1 e^{-B}\cr 
c_2 d_2 e^B & -c_2d_1 e^{-A}\epm,$$
where $$A= (\a_2-\a_1)\xi + (\a_2^{-1}-\a_1^{-1})\eta, \quad B=  
(\a_2+\a_1)\xi +
(\a_2^{-1}+\a_1^{-1})\eta.$$  The new solution $(a, \ti u, \ti v)$ is expressed in terms of $\ti
\pi(\xi,\eta)$ as given in Theorem \ref{ef}. The wave map corresponding to $(a, \ti u, \ti v)$ is
\begin{align*}
\ti s(\xi, \eta)&=\ti E(\xi, \eta, -1)\ti E(\xi, \eta, 1)^{-1}\\
& = h(-1) e^{-ax}\left( \frac{(1+\a_1)(1-\a_2) -
2(\a_1-\a_2) \ti \pi(x,t)}{(1+\a_2)(1-\a_1)}\right) e^{-ax} h(1)^{-1},
\end{align*} 
where $h(\l)= h_{\a_1, \a_2, \pi}(\l)$.  
If $\ti\pi$ has singularities, then
$(a,
\ti u,\ti v)$ and
$\ti s$ have too.  Note that we can choose
$\a_i$ and $v_1, v_2$ so that $h_{\a_1, \a_2, \pi}\bu (a,0,a)$ is singular somewhere or is
smooth on the whole $\R^{1,1}$.  But such wave maps do not have good boundary behavior when
$x\to \pm \infty$.  
\eeg 

\beg\label{ej}
We construct wave maps whose image lies in the $1$-dimensional subgroup
$\R^+=\{\diag(e^t, e^{-t})\n t\in \R\}$ of $SL(2,\R)$.   Let $h(\xi)$ and $k(\eta)$ be smooth real
valued functions, and 
$$a(\xi)= h'(\xi)\diag(1,-1), \quad b(\eta)= k'(\eta) \diag(1, -1).$$
Then $(a, 0, b)$ is a solution of the $-1$ flow equation associated to $SL(2, \R)$, its Lax pair is $\o_\l
= a(\xi)\l\ d\xi + b(\eta)\l^{-1} \ d\eta$, and the trivialization of $\o_\l$ is
$$E(\xi, \eta, \l) = \diag(e^{h(\xi)\l + k(\eta)\l^{-1}}, e^{-(h(\xi)\l + k(\eta)\l^{-1})}).$$
The corresponding wave map is 
\begin{equation}\label{eh}
s(\xi, \eta)= E(\xi, \eta, -1)E(\xi, \eta, 1)^{-1}= \diag(e^{-2(h(\xi) + k(\eta))}, e^{2(h(\xi) +
k(\eta))}).
\end{equation}

Since the subgroup $\R^+$ is abelian, the equation for wave maps into $\R^+$ is the linear wave
equation. Hence every wave maps into $\R^+$ is of the form given in \eqref{eh} for some smooth one
variable real valued functions $h, k$.   
\eeg

\beg 
We compute B\"acklund transformation of the wave maps given in Example \ref{ej}.  
Let $\a_1, \a_2\in \R$, and
$\pi$ the projection of $\C^2$ onto
$\C y_1$ along
$\C y_2$, where $y_1= \bpm c_1\cr d_1\epm$ and
$y_2= \bpm c_2\cr d_2\epm$  are in $\R^2$. 
To get $h_{\a_1, \a_2, \pi}\bu (a, 0, b)$, we first compute 
\begin{align*}
\ti y_1(\xi, \eta)&=E(\xi, \eta, \a_1)^{-1}(y_1) = \bpm c_1 e^{-A_1}\cr d_1 e^{A_1}\epm, \cr
 \ti y_2(\xi, \eta)&= E(\xi, \eta, \a_2)^{-1}(y_2)= \bpm c_2 e^{-A_2}\cr d_2 e^{A_2}\epm,
\end{align*}
where 
$$A_i = h(\xi)\a_i + k(\eta) \a_i^{-1}\quad {\rm for \ }i=1, 2.$$ Let $\ti \pi(\xi, \eta)$ be the 
projection onto $\C \ti y_1(\xi, \eta)$ along $\C \ti y_2(\xi, \eta)$.  Then 
\begin{equation}\label{ei} 
\ti\pi= \frac{1}{W}\bpm c_1d_2 e^{-A_1+A_2}& -c_1c_2 e^{-(A_1+A_2)}\cr 
d_1d_2 e^{A_1+A_2}& - c_2d_1 e^{A_1-A_2}\epm,
\end{equation}   
where 
$$W:= c_1d_2 e^{-A_1+ A_2} - c_2 d_1 e^{A_1-A_2}.$$
So  $(a, \ti u, \ti v):= h_{\a_1, \a_2, \pi}\bu (a, 0, b)$ is given by
$$\bca\ti u= (\a_1-\a_2) [a, \ti\pi], &\cr
\ti v= \left({\rm I\/}-\frac{\a_1-\a_2}{\a_1}\ \ti \pi\right) b \left({\rm I\/}-\frac{\a_1-\a_2}{\a_1}\
\ti \pi\right)^{-1}.\eca$$
The trivialization of $(a, \ti u, \ti v)$ is 
$$\ti E(\xi, \eta, \l) = h_{\a_1, \a_2, \pi}(\l) E(\xi, \eta, \l) h_{\a_1, \a_2, \ti
\pi(\xi,\eta)}^{-1}(\l),$$ and the corresponding wave map is
$$\ti s(\xi, \eta) = h(-1)A_0\left( \frac{(1+\a_1)(1-\a_2) -2(\a_1-\a_2)\ti \pi}{(1+\a_2)(1-\a_1)}\right)
A_0 h(1)^{-1},$$
where $h(\l)= h_{\a_1, \a_2, \pi}(\l)$ and  $A_0(\xi, \eta)= \diag(e^{-(h(\xi) + k(\eta))}, \ e^{(h(\xi) +
k(\eta))})$.

 It can be easily checked that 
\ben 
\item[(i)] $\ti\pi(\xi, \eta)$ has singularity at $(\xi_0, \eta_0)$ if and only if $W(\xi_0,\eta_0)=0$. 
For example,  if $\frac{c_1d_2}{c_2d_1} >0$, then $\ti\pi$ has singularities at points on the
curve
\begin{equation}\label{ek}
h(\xi)(\a_1-\a_2)+ k(\eta)(\a_1^{-1}- \a_2^{-1}) - \frac{1}{2} \ln\left(
\frac{c_1d_2}{c_2d_1}\right)=0.
\end{equation}
\item[(ii)] If $\frac{c_1d_2}{c_2 d_1}<0$, then $W$ never vanishes.  So $\ti\pi$ is smooth for all
$(\xi,
\eta)$.  
 \een

Now suppose that both $h$ and $k$ is in $L^2_1$, i.e, $h, k, h', k'$  are square
integrable.  
If $W(x, t_0)\not=0$ for all $x\in \R$, then the formula for $\ti
\pi$ implies that 
$$\lim_{|x|\to \infty} \ti\pi(x,t_0)= \frac{1}{c_1d_2-c_2d_1}\bpm c_1d_2& -c_1c_2\cr d_1d_2&
-d_1c_2\epm =\pi.$$
So 
$$\lim_{|x|\to -\infty} \ti s(x, t_0)=  h(-1)\left(\frac{(1+\a_1)(1-\a_2) -
2(\a_1-\a_2)\pi}{(1+\a_2)(1-\a_1)}\right) h(1)^{-1}
$$   is a constant

Rewrite the left hand side of \eqref{ek} in space time coordinates $x, t$, and set
$$f(x,t)= (\a_1-\a_2)h((x+t)/2)- (\a_1^{-1}-\a_2^{-1}) k((x-t)/2) - \frac{1}{2} \ln\left(
\frac{c_1d_2}{c_2d_1}\right) $$

For case (i), we can choose $\a_1,\a_2\in \R$ and  $h, k$ in $L^2_1$ so that 
$f(x,0)$ never vanishes  for all $x\in\R$, but vanishes at some $(x_0,t_0)$ for some $t_0>0$.   We
check that the wave map
$\ti s (x,0)$ is smooth with finite energy and constant boundary condition, but it develops singularities in finite
time.  

For case (ii),  we have proved  that $\ti s(x,t)$ has constant boundary condition  for all $t$.  
Claim that the energy of $\ti s$ is finite.  To see this,
note that by Theorem \ref{ah} we have
\begin{align*}
&\tr(\ti s^{-1}\ti s_\xi, \ti s^{-1}s_\xi) = 4\ \tr(a, a), \cr
 &\tr(\ti s^{-1}\ti s_\eta,\ti s^{-1}\ti s_\eta) = 4\ \tr(\ti v, \ti v) 
 \end{align*}
But $\tr(\ti v, \ti v)= \tr(b, b)$, which is finite.  Thus $\ti s$ is smooth, $\lim_{|x|\to
\infty} \ti s(x,t)$ is a fixed constant for all $t$, and  has finite energy. 
 
\eeg

\bs\bs

\end{document}